\documentclass[aap,MSNbibl,dvips]{arximspdf}
\usepackage{mathbh}

\doi{10.1214/13-AAP984} 
\volume{24}
\issue{6}
\pubyear{2014}
\firstpage{2527}
\lastpage{2559}

\makeatletter

\newcommand{\rrvert}{\vert}
\newcommand{\llvert}{\vert}
\newcommand{\eqref}[1]{(\ref{#1})}
\newcommand{\Pd}{\mathbb{P}}

\newcommand{\bbR}{\mathbb{R}}

\newcommand{\bbZ}{\mathbb{Z}}
\newcommand{\bbN}{\mathbb{N}}
\newcommand{\Ex}{\mathbb{E}}

\newcommand{\calO}{\mathcal{O}}
\newcommand{\calA}{\mathcal{A}}

\newcommand{\mB}{\mathcal{B}}
\newcommand{\mC}{\mathcal{C}}

\newcommand{\argmax}{\operatorname{argmax}}

\newcommand{\hY}{\widehat{Y}}
\newcommand{\hS}{\widehat{S}}

\newcommand{\hF}{\widehat{F}}
\newcommand{\hX}{\widehat{X}}
\newcommand{\hE}{\widehat{E}}

\newcommand{\hM}{\widehat{M}}

\newtheorem{thmm}{Theorem}[section]

\newproclaim{defin}{Definition}[section]
\newtheorem{cor}[thmm]{Corollary}
\newproclaim{asum}{Assumption}[section]

\newproclaim{rem}{Remark}[section]
\newtheorem{lem}{Lemma}[section]

\makeatother

\begin{document}
\begin{frontmatter}

\title{Diffusion models and steady-state approximations for
exponentially ergodic Markovian queues}
\runtitle{Steady-state approximations for Markovian queues}

\begin{aug}
\author{\fnms{Itai} \snm{Gurvich}\corref{}\ead[label=e2]{i-gurvich@kellogg.northwestern.edu}}
\runauthor{I. Gurvich}
\address{Kellogg School of Management\\
Northwestern University\\
Evanston, Illinois 60208\\
USA\\
\printead{e2}}
\affiliation{Northwestern University}
\end{aug}

\received{\smonth{4} \syear{2013}}
\revised{\smonth{10} \syear{2013}}

%
\begin{abstract}
Motivated by queues with many servers, we study Brownian steady-state
approximations for continuous time Markov chains (CTMCs). Our
approximations are based on \textit{diffusion models} (rather than a
diffusion limit) whose steady-state, we prove, approximates that of the
Markov chain with notable precision. Strong approximations provide such
``limitless'' approximations for process dynamics. Our focus here is on
steady-state distributions, and the diffusion model that we propose is
tractable relative to strong approximations.

Within an asymptotic framework, in which a scale parameter $n$ is taken
large, a uniform (in the scale parameter) Lyapunov condition imposed on
the sequence of diffusion models guarantees that the gap between the
steady-state moments of the diffusion and those of the properly
centered and scaled CTMCs shrinks at a rate of $\sqrt{n}$.

Our proofs build on gradient estimates for solutions of the Poisson
equations associated with the (sequence of) diffusion models and on
elementary martingale arguments. As a by-product of our analysis, we
explore connections between Lyapunov functions for the fluid model, the
diffusion model and the CTMC.
\end{abstract}

%

\end{frontmatter}

\section{Introduction}\label{sec1}

Fluid and diffusion limits for queuing systems have been applied
successfully toward performance analysis and optimization of various
queuing systems. We are concerned here with performance analysis in
steady-state and, more specifically, with Brownian steady-state
approximations for continuous time Markov chains (CTMCs).

The framework of diffusion limits begins with a sequence of CTMCs $\{
X^n\}$, and properly scaled and centered versions
\[
\hX^n=\frac{X^n-\bar{x}^n}{\sqrt{n}}
\]
for some sequence $\{\bar{x}^n\}$ that arises from the specific
structure of the model. With appropriate assumptions on the parameters
of the CTMC, and on the sequence of initial conditions $\{\hX^n(0)\}$,
one typically proceeds to establish process convergence
%
\begin{equation}
\hX ^n\Rightarrow\hX\qquad \mbox{as } n\rightarrow\infty ,\label{eq:FCLT}
\end{equation}
in the
appropriate function space where $\hX$ is a diffusion process. If each
of the $\{X^n\}$ as well as $\hX$ are ergodic, and $f$ is a continuous
function such that $\{f(\widehat{X}^n(\infty))\}$ is uniformly
integrable, one
can subsequently conclude that
\[
\Ex\bigl[f\bigl(\widehat{X}^n(\infty)\bigr)\bigr]\rightarrow\Ex
\bigl[f\bigl(\hX(\infty )\bigr)\bigr] \qquad\mbox{as } n\rightarrow\infty ,
\]
where $\hX^n(\infty)$ and $\hX(\infty)$ have, respectively, the
steady-state distributions of $\hX^n$ and $\hX$. A relatively general
framework toward proving the required uniform integrability has been
developed in \cite{GZ:06} and applied there to generalized Jackson
networks; see also \cite{BudhirajaLee:07}. It was subsequently applied
successfully to other queueing systems. This so-called \textit{interchange
of limits} establishes that
%
\begin{equation}
\Ex\bigl[f\bigl(\widehat{X}^n(\infty)\bigr)\bigr]=\Ex\bigl[f\bigl(
\hX(\infty )\bigr)\bigr]+o(1)\label{eq:gapintro},
\end{equation}
and supports using $\Ex[f(\hX(\infty))]$ as an approximation for
$\Ex
[f(\widehat{X}^n(\infty))]$.

A central benefit of the limit approach to approximations is the
relative tractability of the diffusion $\hX$ relative to the original
CTMC. The convergence rate embedded in the $o(1)$ term is not, however,
precisely captured by these convergence arguments. In this paper, we
prove that an appropriately defined sequence of diffusion models, that
are as tractable as the diffusion limit, provides accurate
approximations for the steady-state of the CTMCs with an approximation
gap that shrinks at a rate of $\sqrt{n}$. Our approach does not require
process convergence as in \eqref{eq:FCLT}.

We proceed to an informal exposition of the results and key ideas. The
Markov chains that we consider have a semi-martingale representation
\[
X^n(t)=X^n(0)+\int_0^t
F^n\bigl(X^n(s)\bigr)\,ds+M^n(t),
\]
where $M^n$ is a local martingale with respect to a properly defined
filtration. We define a fluid model by (heuristically) removing the
martingale term, that is,
\renewcommand{\theequation}{FM}
\begin{equation}\label{eqFM}
\bar{x}^n(t)=\bar{x}^n(0)+\int
_0^t F^n\bigl(\bar{x}^n(s)
\bigr)\,ds.
\end{equation}
If the FM has a unique stationary point $\bar{x}^n_{\infty}$ satisfying
$F^n(\bar{x}^n_{\infty})=0$, it subsequently makes sense to center
$X^n$ around $\bar{x}_{\infty}^n$ and consider the centered and scaled
process $\widehat{X}^n=(X^n-\bar{x}^n_{\infty})/\sqrt{n}$. The
process $\hX^n$
satisfies the equation
\[
\widehat{X}^n(t)=\widehat{X}^n(0)+ \int
_0^t \hF^n\bigl(\widehat
{X}^n(s)\bigr)\,ds+M^n(t)/\sqrt{n},
\]
where $\hF^n(y)=F^n(\sqrt{n}y+\bar{x}^n_{\infty})/\sqrt{n},   y\in
\bbR
^d$. Under appropriate conditions, a~strong approximation for $\widehat
{X}^n$ is
given by the diffusion process
\[
\widehat{S}^n(t)=\widehat{S}^n(0)+ \int
_0^t \hF^n\bigl(\widehat
{S}^n(s)\bigr)\,ds+\int_0^t
\sigma^n\bigl(\widehat{S}^n(s)\bigr)\,dB(s),
\]
where $B$ is a standard Brownian motion and $\sigma^n$ arises naturally
from the Markov-chain transition functions and is intimately related to
the predictable quadratic variation of the martingale $M^n$.
Strong-approximations theory predicts an approximation gap that is
logarithmic in $nT$ where $T$ is the time horizon; see Remark~\ref{rem:strong}.

A cruder approximation is obtained by replacing the (state dependent)
diffusion coefficient with its value at the stationary point of the FM,
$\bar{x}^n_{\infty}$, to obtain the diffusion process specified by
the equation
\renewcommand{\theequation}{DM}
\begin{equation}\label{eqDM}
\hY^n(t)=\hY^n(0)+ \int_0^t
\hF^n\bigl(\hY^n(s)\bigr)\,ds+\sigma^n\bigl(
\bar {x}^n_{\infty}\bigr)B(t).
\end{equation}
Our main finding is that this straightforward heuristic derivation of
the DM---building on a stationary point of the fluid model to construct
a simplified diffusion model---may provide, insofar as steady-state
analysis is concerned, an impressively accurate approximation.

More precisely, but still proceeding informally at this stage, we prove
the following. Let $\mathcal{A}^n$ be the generator of the diffusion
$\hY^n$. If there exists a function $V$ together with finite positive
constants $b,\delta$ and a compact set $B$ (all not depending on $n$)
such that
\renewcommand{\theequation}{UL}
\begin{equation}\label{eqUL}
\mathcal{A}^nV(x)\leq-\delta V(x)+b\mathbh{1}_B(x),\qquad
x\in\bbR^d,
\end{equation}
then
\[
\Ex\bigl[f\bigl(\hY^n(\infty)\bigr)\bigr]-\Ex\bigl[f\bigl(
\hX^n(\infty)\bigr)\bigr]=\calO(1/\sqrt{n})
\]
for all functions $f$ with $|f|\leq V$. The uniform Lyapunov
requirement UL must be proved on a case-by-case basis, and we
illustrate this via two examples in Section~\ref{sec:examples}. The
requirement UL restricts the scope of our results to (sequences of)
chains in which the corresponding DM is exponentially ergodic.

The sequence of Poisson equations (associated with the sequence of DMs)
is central to our proofs. Let $\pi^n$ be the steady-state distribution
of the diffusion model and $\nu^n$ be that of the scaled CTMC. Let $f$
be such that $\pi^n(f)=0$. [The requirement that $\pi^n(f)=0$ is not
necessary and is imposed in this discussion for expositional purposes.]
We will show that a solution $u^n_f\in\mathcal{C}^{2}(\bbR^d)$ exists
for the DM's Poisson equation
\[
\calA^n u=-f.
\]
Based on It\^{o}'s rule one expects that
\[
\Ex_{\pi^n} \bigl[u^n_f\bigl(
\hY^n(t)\bigr) \bigr]=\Ex_{\pi^n} \bigl[u^n_f
\bigl(\hY ^n(0)\bigr) \bigr]+\Ex_{\pi^n} \biggl[\int
_0^t \calA^n u^n_f
\bigl(\hY ^n(s)\bigr)\,ds \biggr].
\]
Since the DM has, by construction, a diffusion coefficient that does
not depend on the state, the Poisson equation is (for each $n$) a
linear PDE, and we are able to build on existing theory to identify
gradient estimates that are uniform in the index~$n$. These gradient
estimates facilitate proving that
\[
\Ex_{\nu^n} \bigl[u^n_f\bigl(\hX^n(t)
\bigr) \bigr]=\Ex _{\nu^n} \bigl[u^n_f\bigl(
\hX^n(0)\bigr) \bigr]+\Ex_{\nu^n} \biggl[\int
_0^t \calA^n u^n_f
\bigl(\hX ^n(s)\bigr)\,ds \biggr]+t\calO(1/\sqrt{n}).
\]
Informally speaking, this shows that $u_f^n$ ``almost solves'' the
Poisson equation for the CTMC.

Stationarity then allows us to conclude that
\[
\Ex_{\nu^n} \biggl[\int_0^t
\calA^n u^n_f\bigl(\hX^n(s)
\bigr)\,ds \biggr]=-t\Ex_{\nu
^n} \biggl[\int_0^t
f\bigl(\hX^n(s)\bigr)\,ds \biggr]=t\calO(1/\sqrt{n}),
\]
and, in particular, that
\[
\nu^n(f)=\calO(1/\sqrt{n}).
\]
Recalling that $\pi^n(f)=0$, it then follows that
\[
\nu^n(f) -\pi^n(f)=\calO(1/\sqrt{n}).
\]
In the process of proving these results, we explore connections between
the stability of the CTMC and that of the corresponding FM and DM.

Refined properties of the Poisson equation in the context of diffusion
approximations for diffusions with a fast component are used in \cite
{pardoux2001poisson}. In the spirit of this paper, derivative bounds
for certain Dirichlet problems are used in \cite{Excursion} to study
\textit{universal} approximations for the birth-and-death process
underlying the so-called Erlang-A queue. The proofs there are based on
the study of excursions but are closely related to ours; we revisit the
Erlang-A queue in Section~\ref{sec:examples}. The use of gradient
estimates in conjunction with martingale arguments is also the theme
in~\cite{AtaGurvich} where these are used to study optimality gaps in the
control of a multi-class queue. The Poisson equation is replaced there
with the PDE associated with the HJB equation.

\textit{Notation}. Unless stated otherwise, all convergence
statements are for $n\rightarrow\infty $.
We use $|x|$ to denote the Euclidean norm of $x$ in $\bbR^d$ (the
dimension $d$ will be clear form the context). For two nonnegative
sequences $\{a^n\}$ and $\{b^n\}$ we write $a^n=\calO(b^n)$ if
$\limsup_{n\rightarrow\infty }a^n/b^n<\infty$. Throughout we adopt the
convention that
$0/0=0$. We let
\[
B_x(M)=\bigl\{y\in\bbR^d\dvtx |x-y|<M\bigr\},
\]
and denote its closure by $\overline{B}_x(M)$. Following standard
notation, we let $\mC^{j}(\bbR^d)$ be the space of $j$-times
continuously differentiable functions from $\bbR^d$ to $\bbR$, and for
$u\in\mC^2(\bbR^d)$ we let $Du$ and $D^2u$ denote the gradient and the
Hessian of $u$, respectively.

Given a Markov process $\Xi=(\Xi(t), t\geq0)$ on a complete and
separable metric space $\mathcal{X}$, we let $\Pd_x$ be the probability
distribution under which $\Pd\{\Xi(0) = x\} = 1$ for $x\in\mathcal{X}$
and $\Ex_x[\cdot] = \Ex[\cdot|\Xi(0) = x]$ be the expectation operator
w.r.t. the
probability distribution $\Pd_x$. Let $\Pd_{\pi}$ denote the
probability distribution under which $\Xi(0)$ is distributed
according to $\pi$ and put $\Ex_{\pi}[\cdot]$ to be the
expectation operator w.r.t. this distribution. A probability
distribution $\pi$ defined on $\mathcal{X}$ is said to be a
stationary distribution if for every bounded continuous
function $f$
\[
\Ex_{\pi}\bigl[f \bigl(\Xi(t)\bigr)\bigr] = \Ex_{\pi}\bigl[f
\bigl(\Xi(0)\bigr)\bigr]\qquad \mbox{for all } t\geq0.
\]
It is said to be the steady-state distribution if for every such
function and all $x\in\mathcal{X}$,
\[
\Ex_x\bigl[f\bigl(\Xi(t)\bigr)\bigr] \rightarrow\Ex_{\pi}
\bigl[f\bigl(\Xi(0)\bigr)\bigr] \qquad\mbox{as } t\rightarrow\infty.
\]

Given a probability distribution $\nu$ and a nonnegative function $f$,
we define $\nu(f)=\int f(x)\,d\nu(x)$ (which can be infinite). For a
general (not necessarily nonnegative) function, we define $\nu(f)$ as
above whenever $\nu(|f|)<\infty$. Finally, whereas our results
are not concerned with process-convergence, we will be making
connections to the functional central limit theorem. All the processes
that we study are assumed to be right continuous with left limits
(RCLL), and $\Rightarrow$ will be used for convergence in the space
$\mathcal{D}^d[0,\infty)$ of such functions unless otherwise stated.
For RCLL processes we use $x(t-)=\lim_{s\uparrow t}x(s)$ and let
$\Delta x(t)=x(t)-x(t-)$.

\section{A sequence of CTMCs}
We consider a sequence $\{X^n, n\in\bbN\}$ of continuous time Markov
chains (CTMCs). The chain $X^n$ moves on a countable state space
$E^n\subset\mathbb{R}^d$ according to transition rates $\beta
_{y-x}^n(x)=q_{x,y}^n$ for $x,y\in E^n$. Given a nonrandom initial
condition $X^n(0)\in E^n$, the dynamics of $X^n$ are constructed as follows:
\[
X^n(t)=X^n(0)+\sum_{\ell}
\ell Y_{\ell} \biggl(\int_0^t
\beta_{\ell
}^n\bigl(X^n(s)\bigr)\,ds \biggr),
\]
where $\ell\in\mathcal{L}^n=\{y-x\dvtx x,y \in E^n\}$ and $\{Y_{\ell
}, \ell
\in\mathcal{L}^n\}$ are independent unit-rate Poisson processes; see
\cite{EtK:86}, Section~6.4. Letting $\widetilde{Y}_{\ell}(t)=Y_{\ell
}(t)-t$, we rewrite
\[
X^n(t)=X^n(0)+\int_0^t
F^n\bigl(X^n(s)\bigr)\,ds+\sum
_{\ell} \ell\widetilde {Y}_{\ell} \biggl(\int
_0^t \beta_{\ell}^n
\bigl(X^n(s)\bigr)\,ds \biggr),
\]
where
%
\renewcommand{\theequation}{\arabic{equation}}
\setcounter{equation}{2}
\begin{equation}
F_i^n(x)=\sum_{\ell}
\ell_i \beta_{\ell
}^n(x).\label
{eq:Fdefin}
\end{equation}
Provided that $X^n$ is nonexplosive,
\[
M^n(t)=\sum_{\ell}\ell
\widetilde{Y}_{\ell} \biggl(\int_0^t
\beta _{\ell
}^n\bigl(X^n(s)\bigr)\,ds \biggr),
\]
is a local martingale with respect to the filtration
%
\begin{equation}
\quad\mathcal {F}_t^n=\sigma \biggl\{X^n(0),\int
_0^s \beta_{\ell}^n
\bigl(X^n(u)\bigr)\,du, \widetilde{Y}_{\ell} \biggl(\int
_0^s \beta_{\ell}^n
\bigl(X^n(u)\bigr)\,du \biggr);\ell \in\mathcal{L}^n,s\leq t
\biggr\};\label{eq:filtration}
\end{equation}
see \cite{EtK:86}, Theorem~6.4.1.
The local (predictable) quadratic variation of $M^n$ is given by
\[
\bigl\langle M^n \bigr\rangle(t) = \int_0^t
a^n\bigl(X^n(s)\bigr)\,ds,
\]
where
%
\begin{equation}
a_{ij}^n(x)=\sum_{\ell}
\ell_i\ell_j\beta_{\ell
}^n(x).\label
{eq:adefin}
\end{equation}

In essence, $F^n$ and $a^n$ are defined only for values in $E^n$. We
henceforth assume that they are extended to $\bbR^d$ and, with some
abuse of notation, denote by $F^n$ and $a^n$ these extensions. The
requirements that we impose on these extensions will be clear in what follows.

\textit{Fluid models}. Given $x$, we define the $n$th fluid model by
%
\renewcommand{\theequation}{FM}
\begin{equation}
\bar{x}^n(t)=x+\int_0^t
F^n\bigl(\bar{x}^n(s)\bigr)\,ds,
\end{equation}
or, in differential form,
\[
\dot{\bar{x}}^n(t)=F^n\bigl(\bar{x}^n(t)
\bigr),\qquad \bar{x}^n(0)=x.
\]
If $F^n$ is Lipschitz continuous, the fluid model has a solution. We
will assume that there exists a unique $\bar{x}^n_{\infty}$ satisfying
%
\renewcommand{\theequation}{\arabic{equation}}
\setcounter{equation}{5}
\begin{equation}
F^n\bigl(\bar{x}^n_{\infty}\bigr)=0.\label{eq:stationarypoint}
\end{equation}
This
requirement is intimately linked to our Lyapunov requirement; see Lem\-ma~\ref{asum:L}.

\textit{Centered and scaled process.}
Define the processes
%
\begin{equation}
\widehat{X}^n=\frac{X^n-\bar{x}^n_{\infty}}{\sqrt {n}}, \qquad \hM^n=
\frac{M^n}{\sqrt{n}},\label{eq:hxdefin}
\end{equation}
and
denote by $\widehat{E}^n$ the state space of $\widehat{X}^n$. Letting
\[
\hF^n(x)=\frac{F^n(\bar{x}^n_{\infty}+\sqrt{n}x)}{\sqrt{n}},\qquad x\in \bbR^d,
\]
we have
\[
\widehat{X}^n(t)=\widehat{X}^n(0)+\int
_0^t \hF^n\bigl(\widehat
{X}^n(s)\bigr)\,ds +\hM^n(t).
\]
The martingale $\hM^n$ has the local predictable quadratic variation process
%
\begin{equation}
\bigl\langle\hM^n\bigr\rangle(t)=\int_0^t
\bar{a}^n\bigl(\hX ^n(s)\bigr)\,ds,\label
{eq:M_quad_var}
\end{equation}
where
\[
\bar{a}^n(x)=\frac{a^n(\bar{x}^n_{\infty}+x\sqrt{n})}{n},\qquad x\in \bbR^d.
\]

\textit{Assumptions.}
We assume that the jump sizes are bounded uniformly in~$n$:
%
\begin{equation}
\bar {\ell}=\sup_{n}\argmax\bigl\{|\ell|\in
\mathcal{L}^n\bigr\}<\infty,\label
{eq:boundedjumps}
\end{equation}
and that $n$ is sufficiently large so that $\bar
{\ell}/\sqrt{n}\leq1$.

The sequence $\{\hF^n\}$ is assumed to be uniformly Lipschitz, and $\{
\bar{a}^n\}$ is assumed to have linear growth around $0$. Formally,
there exist constants $K_F$, $K_a$ such that, for all $n$,
%
\begin{equation}
\label{eq:Flip} \bigl|\hF^n(x)-\hF^n(y)\bigr|\leq
K_F|x-y|,\qquad x,y\in\bbR^d
\end{equation}
and
%
\begin{equation}
\bigl|\bar{a}^n(x)-\bar{a}^n(0)\bigr|\leq\frac{K_a}{\sqrt {n}}|x|,\qquad x\in
\bbR^d.\label{eq:alip}
\end{equation}
The requirements \eqref{eq:Flip} and $\hF^n(0)=F^n(\bar{x}_n^{\infty
})/\sqrt{n}=0$ guarantee, in particular, that $|\hF^n(x)|\leq
1+K_F|x|$. Condition \eqref{eq:alip} is equivalently stated in terms of
the (unscaled) $a^n$ as
\[
\bigl|a^n(x)-a^n\bigl(\bar{x}^n_{\infty}
\bigr)\bigr|\leq K_a\bigl|x-\bar{x}^n_{\infty
}\bigr|,\qquad x\in
\bbR^d.
\]
We further assume that \emph{$\bar{a}^n(0)$ is positive definite for
each $n$} and that
%
\begin{equation}
\bar{a}^n(0)\rightarrow\bar{a}, \label{eq:quadraticvar_fluid}
\end{equation}
where $\bar{a}$ is itself positive
definite. The matrix $\bar{a}$ is not used in specifying the diffusion
model in Section~\ref{sec:main}, but the assumption of convergence is
used in our proofs, most notably in that of Theorem~\ref
{thmm:expo_ergo}. In various settings, including our own examples in
Section~\ref{sec:examples}, $\bar{a}^n(0)\equiv\bar{a}$ in which case
the convergence requirement is trivially satisfied.

The requirement that the continuous extension $\hF^n$ satisfies the
uniform Lipschitz requirement \eqref{eq:Flip} is a restriction. It
excludes, for example, single-server queueing systems; we revisit this
point in Section~\ref{sec:conclusions}.

\begin{asum} For each $n\in\bbN$, $X^n$ is nonexplosive, irreducible,
positive recurrent and satisfies \eqref{eq:boundedjumps}--\eqref
{eq:quadraticvar_fluid}.
\label{asum:base}
\end{asum}
Positive recurrence and irreducibility imply ergodicity of $X^n$ and,
in particular, the existence of a steady-state distribution (which is
also the unique stationary distribution). In certain cases, positive
recurrence of $X^n$ need not be a priori assumed; see Theorem~\ref
{thmm:CTMC_lyap} and Remark~\ref{rem:AllInOne}.

Assumption~\ref{asum:base} is imposed for the remainder of this paper.

\section{A diffusion model}
\label{sec:main} Recall that $\bar
{x}^n_{\infty}$ is a stationary point for the \textit{fluid model}
%
\renewcommand{\theequation}{FM}
\begin{equation}
\bar{x}^n(t)=\bar{x}^n(0)+\int
_0^t F^n\bigl(\bar{x}^n(s)
\bigr)\,ds,
\end{equation}
and that $\bar{a}^n(0)=a^n(\bar{x}^n_{\infty})/n$. Fix a probability
space and a $d$-dimensional Brownian motion, and let $\hY^n$ be the
strong solution to the SDE
%
\renewcommand{\theequation}{DM}
\begin{equation}
\hY^n(t)=y+\int_0^t
\hF^n\bigl(\hY^n(s)\bigr)\,ds+\sqrt{\bar{a}^n(0)}
B(t).\label{eq:diffusion}
\end{equation}
The existence and uniqueness of a strong
solution follow from the Lipschitz continuity and linear growth of $\hF
^n$ and the constant diffusion coefficient; see, for example, \cite{KaS:91},
Theorems 5.2.5 and 5.2.9.

\begin{rem}[(On strong approximations)]\label
{rem:strong}{The
strong approximation for $\hX^n$ is a diffusion obtained (heuristically
at first) by taking the ``density'' $\bar{a}^n(x)$ of the quadratic
variation in \eqref{eq:M_quad_var} as the diffusion coefficient, to
define the process
\[
\hS^n(t)=y+\int_0^t
\hF^n\bigl(\hS^n(s)\bigr)\,ds+\int_0^t
\sqrt{\bar {a}^n\bigl(\hS ^n(s)\bigr) }\,dB(s).
\]
The process $\hS^n$ provides a ``good'' approximation for the dynamics
of the CTMC in the sense that
\[
\sup_{0\leq t\leq T}\bigl|\widehat{X}^n(t)-\hS^n(t)\bigr|
\leq\Gamma_T^n \log(n),
\]
where $\{\Gamma_T^n\}$ are random variables with exponential tails
(uniformly in $n$); see, for example, \cite{EtK:86}, Chapters 7.5 and
11.3. Given the cruder (state independent) diffusion
coefficient, the DM $\hY^n$ is not likely to be as precise, over finite
horizons, as the strong approximation. In terms of tractability,
however, the analysis of steady-state is simpler for the DM, insofar as
its steady-state distribution (when it exists) involves linear PDEs;
see, for example, \cite{khasminskii2011stochastic}, Chapter~4.9. Our
main result, Theorem~\ref{thmm:main}, shows that this increased
tractability co-exists with an impressive steady-state-approximation
accuracy.}
\end{rem}

\begin{rem}[(On the diffusion model and the diffusion
limit)]{Suppose that, in addition, Assumption~\ref{asum:base}
%
\renewcommand{\theequation}{\arabic{equation}}
\setcounter{equation}{12}
\begin{equation}
\frac{\beta
_{\ell}^n(\bar{x}^n_{\infty}+\sqrt{n}x)-\beta_{\ell}^n(\bar
{x}^n_{\infty
})}{\sqrt{n}}\rightarrow\widehat{\beta}_{\ell}(x),\label
{eq:betaconv}
\end{equation}
uniformly on compact subsets of $\bbR^d$. If $\widehat{X}^n
(0)\Rightarrow y$, then
\[
\widehat{X}^n\Rightarrow\hY,
\]
where $\hY$ is the strong solution to the SDE
\[
\hY(t)=y+\int_0^t \hF\bigl(\hY(s)\bigr)\,ds+
\sqrt{\bar{a}} B(t),
\]
with $\hF(x)=\sum_{\ell}\ell\widehat{\beta}_{\ell}(x)$ and $\bar
{a}$ is
as in \eqref{eq:quadraticvar_fluid}; see \cite{EtK:86}, Theorem~6.5.4.
Given \eqref{eq:betaconv}, requirements (5.9) and (5.14) of that
theorem are trivially satisfied here due to the bounded jumps. The
final requirement in \cite{EtK:86}, Theorem~6.5.4, that $\tau_a=\inf
\{
t\geq0\dvtx |\hY(t)|\geq a\}$ has $\tau_a\rightarrow\infty $ almost
surely, follows
immediately from the fact that $\hY$ is a strong solution. Further, it
is easily proved that $\hY^n\Rightarrow\hY$. Thus, within a
diffusion-limit framework, the DM is consistent with the diffusion
limit in the sense that $\hY^n$ and $\hX^n$ converge to the same limit.}
\end{rem}

For functions $f\in\mC^2(\bbR^d)$, the generator of $\hY^n$ coincides
with the second order differential operator $\mathcal{A}^n$ defined,
for such functions, by
%
\renewcommand{\theequation}{\arabic{equation}}
\setcounter{equation}{13}
\begin{equation}
\mathcal{A}^n f(x) = \sum_{i=1}^d
\hF_i^n(x)\frac
{\partial
}{\partial x_i}f(x)+\frac{1}{2} \sum
_{i,j}^d\bar{a}^n_{ij}(0)
\frac
{\partial^2}{\partial x_i\,\partial x_j}f(x);\label{eq:diff_gen}
\end{equation}
see, for example, \cite{KaS:91}, Proposition~5.4.2.

We next state the uniform Lyapunov assumption. We say that $V\in\mC
^{2}(\bbR^d)$ is a \textit{norm-like} function if $V(x)\rightarrow
\infty $ as
$|x|\rightarrow\infty
$. A function $V\in\mC^2(\bbR^d)$ is said to be sub-exponential if
$V\geq1$ and there exist constants $c_1,c_2$ and $c_3$ such that
%
\begin{equation}
\bigl|DV(x)\bigr|\vee\bigl|D^2V(x)\bigr|\leq c_1e^{c_2|x|},\qquad x\in
\bbR^d\label{eq:expobound1}
\end{equation}
and
%
\begin{equation}
\sup_{y:|y|\leq1}\frac{V(x+y)}{V(x)}\leq c_3,\qquad x\in
\bbR^d.\label{eq:expobound}
\end{equation}

\begin{asum}\label{asum:L} There exist a sub-exponential norm-like function $V \in
\mC
^2(\bbR^d)$ and finite positive constants $b,\delta,K$ (not depending
on $n$) such that
%
\renewcommand{\theequation}{UL}
\begin{equation}
\calA^n V(x)\leq-\delta V(x)+b\mathbh{1}
_{\overline{B}_0(K)}(x) \qquad\mbox{for all } x\in\bbR^d,
\end{equation}
and, for each $n$ and all $x\in\hE^n$,
%
\renewcommand{\theequation}{\arabic{equation}}
\setcounter{equation}{16}
\begin{equation}
\Ex_x \biggl[\int_0^t \bigl(
\bigl(1+\bigl|\widehat{X}^n(s)\bigr|\bigr)^4 V\bigl(
\widehat{X}^n(s)\bigr) \bigr)^2 \,ds \biggr] <\infty, \qquad t
\geq0. \label{eq:finite_integral}
\end{equation}
\end{asum}

Assumption~\ref{asum:L} is imposed for the remainder of this paper. The
requirement that $V\geq1$ is made without loss of generality. If a
norm-like function $V$ satisfies UL, there exists re-defined constants
$b,\delta$ and $K$ such that $1+V$ satisfies UL. All polynomials
$V\geq
1$ satisfy \eqref{eq:expobound1} and \eqref{eq:expobound}---the former
is used only in the proof of Lemma~\ref{lem:fromdisctocont}, and the
latter is used in the derivations of gradient bounds following the
statement of Theorem~\ref{thmm:PDE}. Requirement \eqref
{eq:finite_integral} is relatively unrestrictive as it is imposed on
each individual $n$ (rather than uniformly in $n$).

Lyapunov conditions are frequently used in the context of stability of
continuous time Markov processes (corresponding to fixed $n$ here); see
\cite{TM:93}. The requirement of a uniform Lyapunov condition imposed
on a family of Markov processes is less common (see \cite
{galtchouk2010geometric} for a related example). In Section~\ref
{sec:examples} we study two examples for which all the requirements of
Assumption~\ref{asum:L} are met.

With Assumption~\ref{asum:L}, the existence and uniqueness of a
steady-state distribution, $\pi^n$, for $\hY^n$ follows from \cite{TM:93},
Sections 4 and 6, as does the fact that $\hY^n$ is
exponentially ergodic and that, for each $n$, $\pi^n(|f|)<\infty$ for
all functions $f$ with $|f|\leq V$; see \cite{TM:93}, Theorem~4.2. For
$V$ that satisfies \eqref{eq:expobound1} we have, for all $t\geq0$ and
$x\in\bbR^d$, that
%
\begin{equation}
\Ex_x\bigl[V\bigl(\hY^n(t)\bigr)\bigr]=V(x)+
\Ex_x \biggl[\int_0^{t}
\mathcal{A}^n V\bigl(\hY ^n(s)\bigr)\,ds \biggr];
\label{eq:dynkin}
\end{equation}
see, for example, \cite{klebaner2005introduction}, Theorem~6.3. UL then guarantees that
%
\begin{equation}
\Ex_x\bigl[V\bigl(\hY^n(t)\bigr)\bigr]\leq V(x)+
\Ex_x \biggl[\int_0^t \bigl(-
\delta V\bigl(\hY ^n(s)\bigr)+b\bigr)\,ds \biggr]\label{eq:tbound}
\end{equation}
for all $t\geq0$ and $x\in\bbR^d$
and, consequently, that
%
\begin{equation}
\limsup_{n\rightarrow\infty }\pi^n\bigl(|f|\bigr)\leq\frac
{b}{\delta}\label
{eq:fbound}
\end{equation}
for all functions $f$ with $|f|\leq V$; see also \cite{glynn2006bounding},
Corollary~2.\vadjust{\goodbreak}

Important for our analysis is the following consequence of Assumption~\ref{asum:L}.

\begin{thmm}[(Uniform exponential ergodicity)]\label{thmm:expo_ergo} Let $\pi^n$ be the
steady-state distribution of $\hY^n$. Then there exist finite positive
constants $\mathcal{M}$ and $\mu$ such that
%
\begin{equation}
\sup_{n}\sup_{x\in\bbR^d}\sup_{|f|\leq V}
\frac
{1}{V(x)} \bigl|\Ex _x\bigl[f\bigl(\hY^n(t)\bigr)
\bigr]-\pi^n(f) \bigr|\leq\mathcal{M}e^{-\mu t},\qquad t\geq 0.
\end{equation}
\end{thmm}
Bounds on the convergence rate of exponentially ergodic Markov
processes to their steady-state distribution have been studied
extensively in recent literature. Our proof builds specifically on
\cite
{baxendale2005renewal}. The constants $\mathcal{M}$ and $\mu$ are
related to a minorization condition for the discrete-time process $\{
\hY
^n(m),m\in\bbZ_+\}$. In the standard application, these constants may
depend on $n$. To obtain constants that can be used for all $n\in\bbN$
we must argue that a minorization condition is satisfied uniformly in
$n$; the proof of Theorem~\ref{thmm:expo_ergo} is postponed to
Section~\ref{sec:expo_ergo}.

Theorem~\ref{thmm:expo_ergo} has the following important implication:
fixing a function $f$ with $|f|\leq V$ and $\pi^n(f)=0$, we have for
all $x\in\bbR^d$, that
\[
\sup_{n}\bigl\llvert \Ex_x \bigl[f\bigl(
\hY^n(t)\bigr) \bigr]\bigr\rrvert \leq\mathcal{M} V(x)e^{-\mu t},\qquad
t\geq0,
\]
so that
\[
\sup_{n}\int_0^{\infty}\bigl
\llvert \Ex_x \bigl[f\bigl(\hY^n(s)\bigr) \bigr]\bigr
\rrvert \,ds\leq \mathcal{M} V(x)\int_0^{\infty}e^{-\mu s}
\,ds =C V(x)<\infty
\]
for all $x\in\bbR^d$, where the constant $C$ does not depend on $n$ or
$x$. We conclude that
\[
u_f^n(x)=\int_0^{\infty}
\Ex_x \bigl[f\bigl(\hY^n(s)\bigr) \bigr]\,ds
\]
is a well-defined function of $x\in\bbR^d$ and that, for all $n$,
%
\begin{equation}
\bigl|u_f^n(x)\bigr|\leq CV(x),\qquad x\in\bbR^d.\label{eq:ubound}
\end{equation}
Also, for any
fixed $M>0$ and $n\in\bbN$,
%
\begin{equation}
\sup_{x\in\overline{B}_0(M)}\lim_{t\rightarrow
\infty }\biggl\llvert \int
_0^t \Ex _x\bigl[f\bigl(\hY
^n(s)\bigr)\bigr]\,ds-\int_0^{\infty}
\Ex_x\bigl[f\bigl(\hY^n(s)\bigr)\bigr]\,ds\biggr\rrvert
=0.\label
{eq:Poisson_Convergence}
\end{equation}

Define
%
\begin{equation}
\mathcal{B}_x=B_x \biggl(\frac{1}{1+|x|} \biggr),\qquad x
\in\bbR ^d\label
{eq:mBxdefin}
\end{equation}
and
%
\begin{equation}
\label{eq:barf} \bar{f}(x)= \sup_{y\in\mB_x}\bigl|f(y)\bigr|+ \sup
_{y,z\in\mB_x}\frac{|f(y)-f(z)|}{|y-z|}.
\end{equation}
The introduction of $\bar{f}$ is motivated by the analysis of the
(sequence of) Poisson equations, specifically by the gradient estimates
that require bounds on local fluctuations of $f$; see the derivations
following Theorem~\ref{thmm:PDE}.

Our main result, stated next, establishes that the steady-state
distribution of the Markov chain and the DM are suitably close provided
that moments of the former are uniformly bounded.

\begin{thmm} \label{thmm:main} Fix $V$ that satisfies Assumption~\ref
{asum:L} and a function $f$ such that $\pi^n(f)=0$ and $\bar{f}\leq V$.
Let $\nu^n$ and $\pi^n$ be, respectively, the steady-state
distributions of $\widehat{X}^n$ and $\hY^n$. If
%
\begin{equation}
\limsup_{n\rightarrow\infty } \nu^n \bigl(V(\cdot ) \bigl(1+|
\cdot|\bigr)^4 \bigr)<\infty ,\label{eq:requirement}
\end{equation}
then
\[
\nu^n(f)-\pi^n(f)=\calO (1/\sqrt{n} ).
\]
\end{thmm}

Theorem~\ref{thmm:main} and the remaining results of this section are
proved in Section~\ref{sec:ito}.

\begin{rem}{If $f$ satisfies $\bar{f}\leq V$ but $\pi
^n(f)\neq0$,
consider instead the function $\check{f}^n=f-\pi^n(f)$. Then $\pi
^n(\check{f}^n)=0$. By \eqref{eq:fbound}, $\limsup_{n\rightarrow
\infty }\pi
^n(|f|)\leq b/\delta<\infty$ and, in turn, $\limsup_{n\rightarrow
\infty }\pi
^n(|\check{f}^n|)\leq2b/\delta<\infty$. Further, $\check{f}^n$
satisfies that $\bar{\check{f}}{}^n\leq\bar{f}+\pi^n(|f|)\leq
V+b/\delta
$. Finally, if $V$ satisfies Assumption~\ref{asum:L}, so does the
function $\check{V}=V+b/\delta$. Thus the results that follow hold for
functions $f$ with $\bar{f}\leq V$ regardless of whether $\pi^n(f)=0$
or not.}
\end{rem}

In general, proving requirement \eqref{eq:requirement} (which implies,
in particular, tightness of the sequence $\{\nu^n\}$ of steady-state
distributions) is far from trivial. As we show next \eqref
{eq:requirement} can be argued in advance in our setting. One expects
that, as $n$ grows, the property \eqref{eq:tbound} of the DM will be
approximately valid for the CTMC allowing to draw an implication
similar to \eqref{eq:fbound} with $\hY^n$ there replaced by $\hX^n$.
The next theorem shows that this intuition is valid provided that $V$
satisfies additional simple properties.

Given a function $\Psi\in\mC(\bbR^d)$, define for $x\in\bbR^d$,
%
\begin{equation}
[\Psi]_{2,1,B_x (\bar{\ell}/\sqrt{n}
)}=\sup_{y,z\in
B_x (\bar{\ell}/\sqrt{n} )}\frac{|D^2\Psi(y)-D^2\Psi
(z)|}{|y-z|},\label{eq:square_bracket}\vadjust{\goodbreak}
\end{equation}
where the right-hand side may
be infinite.

\begin{thmm}[{[From DM to CTMC Lyapunov]}]\label{thmm:CTMC_lyap} Let $V$ be as in Assumption~\ref
{asum:L}. Suppose, in addition, that there exists a finite positive
constant $C$ such that, for each $n$, and all $x\in\bbR^d$,
%
\begin{equation}
\bigl(\bigl|DV(x)\bigr|+\bigl|D^2V(x)\bigr|+[V]_{2,1,B_x(\bar{\ell}/\sqrt
{n})}\bigr) \bigl(1+|x|\bigr)\leq
CV(x).\label{eq:DMtoCTMC2}
\end{equation}
Then, for all sufficiently large $n$, and all $x\in\hE^n$,
%
\begin{equation}\quad
\Ex_x\bigl[V\bigl(\widehat{X}^n(t)\bigr)\bigr]\leq V(x)+
\Ex_x \biggl[\int_0^t \biggl(-
\frac
{\delta
}{2}V\bigl(\widehat{X}^n(s)\bigr)+b \biggr) \,ds \biggr],\qquad
t\geq0, \label{eq:lyap_CTMC}
\end{equation}
where $b$ is as in Assumption~\ref{asum:L}. Consequently, $\widehat
{X}^n$ is
ergodic for all such $n$ and, furthermore,
\[
\limsup_{n\rightarrow\infty }\nu^n(V)\leq\frac{2 b}{\delta}.
\]
\end{thmm}

If $V\in\mC^3(\bbR^d)$, condition \eqref{eq:DMtoCTMC2} can be
replaced with
%
\begin{equation}
\bigl(\bigl|DV(x)\bigr|+\bigl|D^2V(x)\bigr|+\bigl|D^3V(x)\bigr|\bigr) \bigl(1+|x|\bigr)\leq
CV(x).\label
{eq:3Diff1}
\end{equation}
Using Taylor's theorem we have, for all $x\in\bbR^d$, that
\begin{eqnarray*}
\bigl(1+|x|\bigr)[V]_{2,1,B_x(\bar{\ell}/\sqrt
{n})}&\leq& \sup_{\eta\in B_x(\bar{\ell}/\sqrt{n})}2\bigl(1+|
\eta|\bigr)\bigl|D^3V(\eta)\bigr|
\\
& \leq& 2 C \Bigl(\sup_{\eta\in B_x(\bar{\ell}/\sqrt{n})}V(\eta) \Bigr)\leq
2c_3CV(x),
\end{eqnarray*}
where the last inequality follows from the sub-exponential property
\eqref{eq:expobound} of $V$ and $\bar{\ell}/\sqrt{n}\leq1$. Note that
\eqref{eq:3Diff1} is satisfied by any polynomial $V\geq1$.

\begin{cor} Fix $V$ that satisfies Assumption~\ref{asum:L}. Suppose
that there exists $\bar{V}$ that, itself, satisfies Assumption~\ref
{asum:L} as well as \eqref{eq:DMtoCTMC2} and
\[
V(\cdot) \bigl(1+|\cdot|\bigr)^4\leq\bar{V}(\cdot).
\]
Then,
\[
\limsup_{n\rightarrow\infty }\nu^n(\bar{V})<\infty,
\]
and, in particular, \eqref{eq:requirement} holds for $V$.\label{cor:tightness}
\end{cor}

\begin{rem}[(A simple case)] {Suppose that $V\in\mC^3(\bbR^d)$ and
satisfies Assumption~\ref{asum:L} and \eqref{eq:3Diff1}. If there
exists $m\in\bbN$ such that $V_m(\cdot)=(V(\cdot))^m\geq V(\cdot
)(1+|\cdot|)^4$ and $V_m$ satisfies \eqref{eq:finite_integral}, then we
can take $\bar{V}=V_m$ in Corollary~\ref{cor:tightness}. Indeed, for an
integer $m\geq2$,
\begin{eqnarray*}
\mathcal{A}^n V_m(x) &=& m V_{m-1}(x)
\mathcal {A}^n \Psi(x) +m(m-1)V_{m-2}(x)\frac{1}{2} \sum_{i,j}^d
\bar {a}^n_{ij}(0)\frac{\partial}{\partial x_i} V(x)
\\
&\leq&-\delta m V_m(x)+bmV_{m-1}(x)+m(m-1)CV_{m-1}(x),
\end{eqnarray*}
with $\delta$ and $b$ as in Assumption~\ref{asum:L} and $C$ as in
\eqref
{eq:3Diff1}. Thus if $V\in\mC^3(\bbR^d)$ is sub-exponential and
satisfies UL and \eqref{eq:3Diff1}, so does $V_m$. \label
{rem:simplecase}}
\end{rem}

\begin{rem}[(A unified set of conditions)]Combined, Theorem~\ref
{thmm:main} and Corollary~\ref{cor:tightness} establish the following:
If there exist functions $V$ and $\bar{V}$ both satisfying Assumption~\ref{asum:L} such that \eqref{eq:DMtoCTMC2} holds for $\bar{V}$ and
$V(\cdot)(1+|\cdot|)^4\leq\bar{V}(\cdot)$, then we \emph{simultaneously} have: (i) the positive recurrence of $\hX^n$ for
sufficiently large $n$, (ii) the moment bound in \eqref{eq:requirement}
(which implies, in particular, the tightness of $\nu^n$) and (iii) the
$\calO(1/\sqrt{n})$ approximation gap.

With the exception of the simple requirement \eqref
{eq:finite_integral}, this reduces the requirements to properties of
the DM. \label{rem:AllInOne}
\end{rem}

We conclude this section with an observation pertaining to the
connection between the stability of the FM and the DM. Suppose that
there exist a norm-like function $V$ and a constant $\eta$ such that
%
\begin{equation}
\label{eq:FMLyap}V (x )>V(0)\quad\mbox{and}\quad \hF ^n(x)'DV(x)
\leq-\eta\bigl(V(x)-V(0)\bigr),\qquad x\neq0.
\end{equation}
Letting $V^n(x)=V
(\frac{x-\bar{x}^n_{\infty}}{\sqrt{n}} )-V(0)$ we have
\[
F^n(x)'DV^n(x)\leq-\eta
V^n(x),\qquad x\neq\bar{x}_{\infty}^n,
\]
so that the FM is stable in the sense that, for each $n$ and any
initial condition $\bar{x}^n(0)\in\bbR^d$,
$\bar{x}^n(t)\rightarrow\bar{x}^n_{\infty}$ as $t \rightarrow
\infty $. Moreover,
%
\begin{eqnarray}
\mathcal{A}^n V(y)& \leq&\hF^n(y)'DV(y)+\bigl|
\bar {a}^n(0)\bigr|\bigl|D^2V(y)\bigr|
\nonumber
\\[-8pt]
\\[-8pt]
\nonumber
& \leq&-\eta\bigl(V(y)-V(0)\bigr)+\bigl|\bar {a}^n(0)\bigr|\bigl|D^2V(y)\bigr|.
\end{eqnarray}
The following is an immediate consequence.
%
\begin{lem}[{[FM and DM stability]}] Let $V\in\mC^2(\bbR^d)$ be a
sub-exponential norm-like function satisfying \eqref
{eq:finite_integral} and \eqref{eq:FMLyap}. If
\[
\limsup_{|x|\rightarrow\infty }\frac{|D^2V(x)|}{V(x)}=0,
\]
then $V$ satisfies UL and, in turn, Assumption~\ref{asum:L}. \label
{lem:FMtoDM}
\end{lem}

\section{A sequence of Poisson equations}

In what follows, fixing a set $\mB\subseteq\mathbb{R}^d$, $\mC
^2(\mB)$
denotes the space of twice continuously differentiable functions from
$\mB$ to $\mathbb{R}$. For $u\in\mC^{2}(\mB)$, recall that $Du$ and
$D^2u$ denote the gradient and the Hessian of $u$, respectively. The
space $\mC^{2,1}(\mB)$ is then the subspace of $\mC^{2}(\mB)$ members
of which have second derivatives that are Lipschitz continuous on $\mB
$. That is, a twice continuously differentiable function $u\dvtx \mathbb
{R}^d\to\mathbb{R}$ is in $\mC^{2,1}(\mB)$ if
\[
[u]_{2,1,\mB}=\sup_{x,y\in\mB, x\neq y} \frac
{|D^2u(x)-D^2u(y)|}{|x-y|}<\infty.
\]
[In equation \eqref{eq:square_bracket} the set $\mB$ is taken to be
$B_x(\bar{\ell}/\sqrt{n})$.] We define
$d_x=\operatorname{dist}(x,\partial\mB)=\inf\{|x-y|, y\in\partial\mB\}$ where
$\partial\mB$ stands for the boundary of $\mB$, and we let
$d_{x,z}=\min
\{d_x,d_z\}$. We define
%
\begin{equation}
|u|_{2,1,\mB}^*= \sum_{j=0}^2
[u]_{j,\mB
}^*+\sup_{x,y\in\mB,x\neq y}\,d_{x,y}^{3}
\frac{|D^2u(x)-D^2
u(y)|}{|x-y|} ,\label{eq:ustardefin}
\end{equation}
where $[u]_{j,\mB}^*=\sup_{x\in
\mB}d_x^j |D^j u(x)|$ for $j=0,1,2$. Above $d_x^j$ (resp., $d_{x,y}^j$)
denotes the $j$th power of $d_x$ (resp., of $d_{x,y}$). We let
$|u|_{0,\mB}=[u]_{0,\mB}^*=\sup_{x\in\mB}|u(x)|$, and
\[
|f|_{0,1,\mB}^{(2)}= \sup_{x\in\mB}
d_x^2 \bigl|f(x)\bigr|+\sup_{x,y\in\mB}d_{x,y}^3
\frac{|f(x)-f(y)|}{|x-y|}.
\]
We say that the function is locally Lipschitz if
$|f|_{0,1,\mB_x}^{(2)}<\infty$ for all $x\in\bbR^d$, where $\mB_x$ is
as in \eqref{eq:mBxdefin}.
%
\begin{thmm}\label
{thmm:PDE} Fix $V$ that satisfies Assumption~\ref{asum:L} and a
locally Lipschitz function $f$ with $|f|\leq V$ and $\pi^n(f)=0$. Then,
for each $n$, the Poisson equation
%
\begin{equation}
\calA^n u=-f\label
{eq:Poisson}
\end{equation}
has a unique solution $u_f^n\in\mC^{2}(\bbR^d)$ given
by
%
\begin{equation}
u_f^n(x)=\int_0^{\infty}
\Ex_x\bigl[f(\hY^n(t)\bigr]\,dt.\label
{eq:poissonsolution}
\end{equation}
Moreover, there exist a finite positive constant $\Theta$ (not
depending on $n$) such that
\[
\bigl|u_f^n\bigr|_{2,1,\mathcal{B}_x}^*\leq\Theta
\bigl(\bigl|u_f^n\bigr|_{0,\mB
_x}+|f|^{(2)}_{0,1,\mathcal{B}_x}
\bigr),\qquad x\in\bbR^d.
\]
Consequently, for all $n$ and $x\in\bbR^d$,
%
\begin{eqnarray}
\bigl|Du_f^n(x)\bigr|&\leq&2\Theta\bigl(\bigl|u_f^n\bigr|_{0,\mB
_x}+|f|^{(2)}_{0,1,\mathcal
{B}_x}
\bigr)\bigl (1+|x|\bigr),\label{eq:PDE0}
\\
\bigl|D^2u_f^n(x)\bigr|&\leq& 4\Theta
\bigl(\bigl|u_f^n\bigr|_{0,\mB
_x}+|f|^{(2)}_{0,1,\mathcal{B}_x}
\bigr) \bigl(1+|x|\bigr)^2\label{eq:PDE1}
\end{eqnarray}
and
%
\begin{equation}
\bigl[u_f^n\bigr]_{2,1,\mB_x}\leq8\Theta
\bigl(\bigl|u_f^n\bigr|_{0,\mB
_x}+|f|^{(2)}_{0,1,\mathcal{B}_x}
\bigr)\bigl (1+|x|\bigr)^3.\label{eq:PDE2}
\end{equation}
\end{thmm}

Several observations are useful for what follows: recall \eqref
{eq:ubound} that $|u_f^n(x)|\leq CV(x)$ for some constant $C$. By the
assumed sub-exponentiality of $V$
\[
\bigl|u_f^n\bigr|_{0,\mB_y}\leq\sup_{z\in\mB_{y}}C
V(z)\leq c_3 C V(y)
\]
for all $y\in\bbR^d$, where $c_3$ is as in \eqref{eq:expobound}. In turn,
\[
\sup_{y\in B_x (\bar{\ell}/\sqrt{n} )}\bigl|u_f^n\bigr|_{0,\mB
_y}
\leq\sup_{y\in B_x (\bar{\ell}/\sqrt{n} )} c_3 CV(y)\leq c_3^2C
V(x).
\]
For a function $f$ with $\bar{f}\leq V$ [see \eqref{eq:barf}] and for
all $y\in\bbR^d$,
\[
|f|_{0,1,\mB_y}^{(2)}\leq\bar{f}(y)\leq V(y),
\]
so that
\[
\sup_{y\in B_x(\bar{\ell}/\sqrt{n})}|f|_{0,1,\mB_y}^{(2)}\leq
c_3 V(x)
\]
for all $x\in\bbR^d$. Defining
%
\begin{equation}
C_V(x)=16\Theta\bigl(1+c_3^2 C\bigr) V(x)
\bigl(1+|x|\bigr)^3,\qquad x\in \bbR^d,\label{eq:Cfdefin}
\end{equation}
we have, by Theorem~\ref{thmm:PDE} (and
assuming, without loss of generality that $c_3\geq1$), that for all
$n\in\bbN$ and $x\in\bbR^d$,
%
\begin{eqnarray}
\label{eq:gradient_est} \bigl |Du_f^n(x)\bigr|&\leq& C_V(x)/\bigl(1+|x|\bigr)^2,
\nonumber\\
\bigl|D^2u_f^n(x)\bigr|&\leq& C_V(x)/\bigl(1+|x|\bigr)\quad
\mbox{and}
\\
 \bigl[u_f^n\bigr]_{2,1,B_x(\bar{\ell}/\sqrt
{n})}&\leq&
C_V(x).\nonumber
\end{eqnarray}

\begin{pf*}{Proof of Theorem \protect\ref{thmm:PDE}} We first prove
that $u_f^n$
in \eqref{eq:poissonsolution} solves the Poisson equation \eqref
{eq:Poisson}. Since $f$ is fixed throughout we omit it from the notation.

Fixing $M$, let $u_M^n$ be the solution to Dirichlet problem
\begin{eqnarray*}
\calA^n u(x) &=&-f(x), \qquad x\in B_0(M);
\\
u &=& u^n,\qquad x\in\partial B_0(M).
\end{eqnarray*}
In the boundary condition, $u^n$ is
as in \eqref{eq:poissonsolution}. The existence and uniqueness of a
solution $u_M^n\in\mC^{0}(\overline{B}_0(M))\cap\mC^{2,1}(B_0(M))$
follows directly from \cite{TandG}, Theorem~6.13, recalling that $\hF^n$
is Lipschitz continuous and $\bar{a}^n(0)$ is a constant matrix and
hence trivially Lipschitz. Theorem~6.13 of \cite{TandG} requires that
$u_M^n$ is continuous in $x$ on $\partial B_0(M)$. This follows exactly
as in part (c) of \cite{pardoux2001poisson}, Theorem~1, using \eqref
{eq:Poisson_Convergence}. We omit the detailed argument.

It follows that
\[
u_M^n(x)=\Ex_x \biggl[\int
_0^{\tau_M^n}f\bigl(\hY^n(s)\bigr)\,ds
\biggr],
\]
where $\tau_M^n=\inf\{t\geq0\dvtx \hY^n(t)\notin B_0(M)\}$; see \cite{KaS:91},
Proposition~5.7.2 and Lem\-ma~5.7.4. We have that
\[
u_M^n(x)=u^n(x)\qquad\mbox{for all } x\in
B_0(M),
\]
with $u^n(x)$ as in
\eqref{eq:poissonsolution}. This assertion is proved as in \cite{pardoux2001poisson},
Theorem~1, part (d). Since $M$ is arbitrary we
conclude that, $u^n(x)$ solves the Poisson equation~\eqref{eq:Poisson}.

To establish the gradient estimates observe that, since $\bar{a}^n(0)$
is bounded in $n$, there exists a constant $C_a$ (not depending on $n$)
such that (with the notation in~\cite{TandG}, Theorem~6.2) $|\bar
{a}^n(0)|_{0,1,\mB_x}^{(0)}\leq C_a$. From the positive definiteness of
$\bar{a}^n(0)$, and since $\bar{a}^n(0)\rightarrow\bar{a}$ for a
positive definite $\bar{a}$, it follows that there exists a constant
$\lambda>0$ such that
%
\begin{equation}
\sum_{ij}\bar{a}_{ij}^n(0)
\xi_i\xi_j\geq \lambda|\xi|^2\label{eq:lambda}
\end{equation}
for all $n$ and all $\xi\in
\bbR^d$.
Finally, following the notation in \cite{TandG}, Theorem~6.2,
\begin{eqnarray*}
\bigl|\hF^n\bigr|_{0,1,\mB_x}^{(1)}&=& \bigl|\hF^n\bigr|_{0,\mB_x}^{(1)}+
\bigl[\hF ^n\bigr]_{0,1,\mB
_x}^{(1)}
\\
&=& \bigl[\hF^n\bigr]_{0,\mB_x}^{(1)}+\sup
_{y,z\in\mB_x} d_{y,z}^2 \frac
{|\hF
^n(y)-\hF^n(z)|}{|y-z|}
\\
&=& \sup_{y\in\mB_x} d_y\bigl|\hF^n(y)\bigr|+\sup
_{y,z} d_{y,z}^2 \frac
{|\hF
^n(y)-\hF^n(z)|}{|y-z|}
\\
& \leq& 2 K_F,
\end{eqnarray*}
where $K_F$ is as in \eqref{eq:Flip}. In turn, by \cite{TandG}, Theorem~6.2, that
\[
\bigl|u_f^n\bigr|_{2,1,\mathcal{B}_x}^*\leq\Theta
\bigl(\bigl|u_f^n\bigr|_{0,\mB
_x}+|f|^{(2)}_{0,1,\mathcal{B}_x}
\bigr),
\]
where $\Theta$ depends only on $K_F,C_a,d$ and the constant $\lambda$
in \eqref{eq:lambda} (for $\Lambda$ there, we take $K_F\vee C_a$).
Bounds \eqref{eq:PDE0}--\eqref{eq:PDE2} now follow from the definition
of $|u_f^n|_{2,1,\mathcal{B}_x}^*$ applied to points in the subset
$B_x(1/(2(1+|x|)))$ of $\mB_x$.
Specifically, for each $y\in\mathcal{B}_x$,
\[
d_y\bigl |D u_f^n(y)\bigr|\leq[u]_{1,\mB_x}^*
\leq\bigl|u_f^n\bigr|_{2,1,\mathcal{B}_x}^*.
\]
Noting that $d_y\geq1/(2(1+|x|))$ for all $y\in B_x(1/(2(1+|x|)))$ we
have, for all such~$y$ (in particular for $x$ itself), that
\[
\bigl|D u_f^n(y)\bigr|\leq\bigl|u_f^n\bigr|_{2,1,\mathcal{B}_x}^*\bigl(1+|x|\bigr).
\]
Equations \eqref{eq:PDE1} and \eqref{eq:PDE2} are argued similarly.
\end{pf*}

\section{Proofs of Theorems \texorpdfstring{\protect\ref{thmm:main}}{3.2} and \texorpdfstring{\protect\ref{thmm:CTMC_lyap}}{3.3}}
\label{sec:ito}

The following simple lemma is proved in the \hyperref[app]{Appendix}. Given a function
$\Psi\in\mC^{2}(\bbR^d)$ we write $\Psi_i$ for the $i$th
coordinate of
$D\Psi$ and $\Psi_{ij}$ for the $ij$th coordinate of $D^2\Psi$.
%
\begin{lem}\label{lem:ito} Let $\Psi\in\mC^{2}(\bbR^d)$ be such that, for all
$x\in
\widehat{E}^n$ and $t\geq0$,
%
\begin{eqnarray}\label{eq:integrability_cond}
&&\Ex_x \biggl[\int_0^t \bigl(\bigl|D
\Psi\bigl(\widehat {X}^n(s)\bigr)\bigr|+\bigl|D^2\Psi\bigl(
\widehat{X}^n (s)\bigr)\bigr|
\nonumber
\\[-8pt]
\\[-8pt]
\nonumber
&&\hspace*{80pt}{}+
[\Psi ]_{2,1,B_{\widehat{X}^n(s)}(\bar{\ell}/\sqrt{n})} \bigr) \bigl(1+\bigl |
\widehat{X}^n (s)\bigr| \bigr) \,ds \biggr]<\infty.
\end{eqnarray}
Then, for all $x\in\widehat{E}^n$ and $t\geq0$,
%
\begin{equation}\qquad
\Ex_x \bigl[\Psi\bigl(\widehat{X}^n(t)\bigr) \bigr]=
\Psi(x)+\Ex_x \biggl[\int_0^t
\calA^n \Psi \bigl(\widehat{X}^n(s)\bigr)\,ds
\biggr]+A_{\Psi}^{n,x}(t)+D_{\Psi
}^{n,x}(t),\label{eq:almostito}
\end{equation}
where $\calA^n$ is as in \eqref{eq:diff_gen} and, for all $x\in
\widehat
{E}^n$ and $t\geq0$,
\begin{eqnarray*}
\bigl|A_{\Psi}^{n,x}(t)\bigr|&\leq& \frac{\bar{\ell}}{2\sqrt{n}}\Ex_x
\biggl[ \int_0^t [\Psi]_{2,1,B_{\widehat{X}^n(s)}(\bar{\ell}/\sqrt{n})}\bigl|\bar
{a}^n\bigl(\widehat{X}^n(s)\bigr)\bigr| \,ds \biggr],
\\
D_{\Psi}^{n,x}(t)&=&\frac{1}{2}\Ex_x
\Biggl[\sum_{i,j}^d \int
_0^t \Psi _{ij}\bigl(
\widehat{X}^n(s)\bigr) \bigl(\bar{a}_{ij}^n
\bigl(\widehat{X}^n(s)\bigr)-\bar {a}_{ij}^n(0)
\bigr)\,ds \Biggr].
\end{eqnarray*}
\end{lem}

Below $\bar{f}$ is as in \eqref{eq:barf} and $C_V$ as in \eqref{eq:Cfdefin}.

\begin{cor}\label{cor:interim} Fix $V$ that satisfies Assumption~\ref{asum:L} and a
function $f$ such that $\bar{f}\leq V$. Then there exists a finite
positive constant $C$ (not depending on $n$), such that, for all $x\in
\widehat{E}^n$ and $t\geq0$,
\begin{eqnarray*}
&&\biggl|\Ex_x \bigl[u_f^n\bigl(
\widehat{X}^n(t)\bigr) \bigr]- u_f^n(x)-\Ex
_x \biggl[ \int_0^t
\calA^n u_f^n\bigl(\widehat{X}^n(s)
\bigr)\,ds \biggr] \biggr|
\\
&&\qquad\leq C \biggl(\Ex _x \biggl[\int_0^t
\frac{C_V(\widehat{X}^n(s))}{\sqrt{n}} \biggl(1+\frac{|\widehat
{X}^n(s)|}{\sqrt {n}} \biggr)\,ds \biggr] \biggr).
\end{eqnarray*}
\end{cor}

\begin{pf}By \eqref{eq:gradient_est} we have, for
$x\in\bbR^d$, that
\begin{eqnarray*}
\bigl(\bigl|Du^n_f(x)\bigr|+\bigl|D^2u^n_f(x)\bigr|+
\bigl[u^n_f\bigr]_{2,1,B_x
({\bar{\ell}}/{\sqrt{n}} )}\bigr) \bigl(1+|x|\bigr) & \leq&3
C_V(x) \bigl(1+|x|\bigr)
\\
& \leq&\varepsilon\bigl(1+|x|\bigr)^4V(x)
\end{eqnarray*}
for some finite positive constant. By Assumption~\ref{asum:L},
specifically \eqref{eq:finite_integral},
\[
\Ex_x \biggl[\int_0^t\bigl(1+\bigl|
\widehat{X}^n(s)\bigr|\bigr)^4 \bigl(V\bigl(\widehat
{X}^n(s)\bigr)\bigr)^2 \,ds \biggr]<\infty,
\]
so that $V$ satisfies the requirements of Lemma~\ref{lem:ito}, and we
have that
%
\begin{eqnarray}
\label{eq:Dbound} \bigl|D_{u^n_f}^{n,x}(t)\bigr|& \leq&\frac
{1}{2}\Ex
_x \biggl[\int_0^t
\bigl|D^2 u^n_f\bigl(\widehat{X}^n(s)
\bigr)\bigr|\bigl|\bar{a}^n\bigl(\widehat {X}^n(s)\bigr)-\bar
{a}^n(0)\bigr|\,ds \biggr]\nonumber
\\
& \leq&\frac{K_a}{2\sqrt{n}}\Ex _x \biggl[ \int
_0^t\bigl |D^2 u^n_f
\bigl(\widehat{X}^n(s)\bigr)\bigr|\bigl|\widehat{X}^n(s)\bigr|\,ds \biggr]
\\
& \leq& \frac{K_a}{2\sqrt{n}}\Ex_x \biggl[\int_0^t
C_V\bigl(\widehat {X}^n(s)\bigr)\,ds \biggr].
\nonumber
\end{eqnarray}
The second inequality follows from \eqref{eq:alip}. 
The last inequality follows from~\eqref{eq:gradient_est}.

Next,
\begin{eqnarray}
\label{eq:Abound} \bigl|A_{u^n_f}^{n,x}(t)\bigr|&\leq&\frac
{\bar
{\ell}}{2\sqrt{n}}
\Ex_x \biggl[\int_0^t
\bigl[u^n_f\bigr]_{2,1,B_{\widehat
{X}^n(s)}(\bar
{\ell
}/\sqrt{n})}\bigl|\bar{a}^n
\bigl(\widehat{X}^n(s)\bigr)\bigr|\,ds \biggr]
\nonumber\\
& \leq &\frac{\bar{\ell}}{2\sqrt{n}}\Ex_x \biggl[\int
_0^t \bigl[u^n_f
\bigr]_{2,1,B_{\widehat{X}^n
(s)}(\bar{\ell}/\sqrt{n})}\bigl|\bar{a}^n(0)\bigr|\,ds \biggr]
\\
&&{} + \frac{\bar{\ell}}{2\sqrt{n}}\Ex_x \biggl[\int_0^t
\bigl[u^n_f\bigr]_{2,1,B_{\widehat{X}^n(s)}(\bar{\ell}/\sqrt{n})}\bigl|\bar {a}^n
\bigl(\widehat{X}^n (s)\bigr)-\bar {a}^n(0)\bigr|\,ds \biggr].
\nonumber
\end{eqnarray}
Using \eqref{eq:alip}, \eqref{eq:quadraticvar_fluid} and \eqref
{eq:gradient_est} we conclude that
\[
\bigl|A_{u^n_f}^{n,x}(t)\bigr|\leq\frac{\bar{\ell}}{2\sqrt {n}}\Ex _x
\biggl[\int_0^t C_V\bigl(
\widehat{X}^n(s)\bigr) \bigl(\bigl|\bar{a}^n(0)\bigr|+K_a\bigl|
\widehat {X}^n(s)\bigr|/\sqrt {n}\bigr)\,ds \biggr],
\]
which completes the proof.
\end{pf}

We are ready to prove Theorem~\ref{thmm:main}.

\begin{pf*}{Proof of Theorem \protect\ref{thmm:main}}
As $\nu^n$ is a stationary distribution we have, by \eqref{eq:ubound}
and \eqref{eq:requirement}, that
\[
\Ex_{\nu^n}\bigl[u_f^n\bigl(
\widehat{X}^n(t)\bigr)\bigr]=\Ex_{\nu^n}\bigl[u_f^n
\bigl(\widehat {X}^n(0)\bigr)\bigr]\leq C\nu ^n(V)<\infty
\]
for all sufficiently large $n$ and all $t\geq0$. Recalling that
$\calA^n u_f^n=-f$, Corollary~\ref{cor:interim} guarantees the
existence of a finite positive constant $\vartheta$ (not depending on
$n$) such that
%
\begin{eqnarray}\label{eq:inter1}\quad
\biggl\llvert \Ex_{\nu^n} \biggl[\int_0^t
f\bigl(\widehat {X}^n(s)\bigr)\,ds \biggr] \biggr\rrvert &\leq&\vartheta
\Ex_{\nu^n} \biggl[\int_0^t
\frac{C_V(\widehat{X}^n
(s))}{\sqrt {n}} \biggl(1+\frac{|\widehat{X}^n(s)|}{\sqrt{n}} \biggr)\,ds \biggr]
\nonumber
\\[-8pt]
\\[-8pt]
\nonumber
&=& \vartheta t\Ex_{\nu^n} \biggl[\frac
{C_V(\widehat{X}^n(0))}{\sqrt{n}} \biggl(1+
\frac{|\widehat
{X}^n(0)|}{\sqrt{n}} \biggr) \biggr]
\end{eqnarray}
for all $t\geq0$, where the interchange of integral and expectation is
justified by the nonnegativity of the integrands. Using again \eqref
{eq:requirement} and the nonnegativity of $V$ we have, for all $t\geq
0$, that
\[
\Ex_{\nu^n} \biggl[\int_0^t\bigl |f\bigl(
\widehat{X}^n(s)\bigr)\bigr|\,ds \biggr]\leq\Ex _{\nu^n} \biggl[ \int
_0^t V\bigl(\widehat{X}^n(s)
\bigr)\,ds \biggr]= t\nu^n(V)<\infty.
\]
This justifies replacing integral and expectation in \eqref{eq:inter1}
to conclude that, with $t>0$,
\begin{eqnarray*}
\bigl|\nu^n(f)\bigr|&=&\frac{1}{t}\biggl\llvert \Ex_{\nu^n}
\biggl[\int_0^t f\bigl(\widehat{X}^n(s)
\bigr)\,ds \biggr]\biggr\rrvert  \leq\vartheta\Ex_{\nu^n} \biggl[
\frac
{C_V(\widehat{X}^n
(0))}{\sqrt{n}} \biggl(1+\frac{|\widehat{X}^n(0)|}{\sqrt{n}} \biggr) \biggr]
\\
&=&\calO (1/\sqrt{n})
\end{eqnarray*}
for a (re-defined) constant $\vartheta$ as required, where the last
equality follows from \eqref{eq:requirement} recalling the definition
of $C_V$ in \eqref{eq:Cfdefin}.
\end{pf*}

\begin{pf*}{Proof of Theorem \protect\ref{thmm:CTMC_lyap}} Let $V$ be
as in
Assumption~\ref{asum:L}. Applying Lemma~\ref{lem:ito} as in the proof
of Corollary~\ref{cor:interim} we have that
\begin{eqnarray*}
\bigl|A_V^{n,x}(t)\bigr|&\leq& \frac{\bar{\ell}}{2\sqrt{n}}\Ex_x
\biggl[\int_0^t [V]_{2,1,B_{\widehat{X}^n
(s)}
({\bar{\ell}}/{\sqrt{n}} )}\bigl|
\bar{a}^n(0)\bigr|\,ds \biggr]
\\
&&{} + \frac{\bar{\ell}}{2\sqrt{n}}\Ex_x \biggl[\int_0^t
[V]_{2,1,B_{\widehat{X}^n(s)} ({\bar{\ell}}/{\sqrt
{n}}
)}\bigl|\bar {a}^n\bigl(\widehat{X}^n(s)
\bigr)-\bar{a}^n(0)\bigr|\,ds \biggr]
\\
& \leq&\Ex_x \biggl[\int_0^t
\frac{\delta}{4} V\bigl(\widehat {X}^n(s)\bigr)\,ds \biggr]
\end{eqnarray*}
for all sufficiently large $n$. The last inequality follows noting
that, by \eqref{eq:alip}, \eqref{eq:quadraticvar_fluid} and~\eqref{eq:DMtoCTMC2}, there exists a finite positive constant $C$ such that
\[
[V]_{2,1,B_{\widehat{X}^n(s)} ({\bar{\ell}}/{\sqrt
{n}}
)}\bigl|\bar {a}^n(0)\bigr|\leq C V\bigl(
\widehat{X}^n(s)\bigr)
\]
and
\[
[V]_{2,1,B_{\widehat{X}^n(s)} ({\bar{\ell}}/{\sqrt
{n}}
)}\bigl|\bar {a}^n\bigl(\widehat{X}^n(s)
\bigr)-\bar{a}^n(0)\bigr|\leq\frac{CK_a}{\sqrt {n}}V\bigl(\widehat{X}^n(s)
\bigr),
\]
where $K_a$ is as in \eqref{eq:alip}. Similarly one argues, using
\eqref
{eq:alip} and \eqref{eq:DMtoCTMC2}, that for all sufficiently large $n$,
\begin{eqnarray*}
\bigl|D_V^{n,x}(t)\bigr| & \leq&\frac{1}{2}\Ex_x
\biggl[\int_0^t \bigl|D^2V\bigl(
\hX^n(s)\bigr)\bigr|\bigl|\bar{a}^n\bigl(\hX^n(s)\bigr)-
\bar{a}^n(0)\bigr|\,ds \biggr]
\\
& \leq&\Ex_x \biggl[\int_0^t
\frac{\delta}{4} V\bigl(\widehat {X}^n(s)\bigr)\,ds \biggr],
\end{eqnarray*}
to conclude from Assumption~\ref{asum:L} and Lemma~\ref{lem:ito} that
\[
\Ex_x\bigl[V\bigl(\widehat{X}^n(t)\bigr)\bigr]\leq
V(x)+\Ex_x \biggl[\int_0^t
\biggl(-\frac{\delta
}{2}V\bigl(\widehat{X}^n(s)\bigr)+b \biggr)\,ds
\biggr].
\]
In turn, \eqref{eq:lyap_CTMC} holds for all sufficiently large $n$.

This guarantees that $\hX^n$ is ergodic for all such $n$; see, for
example, \cite{robert2003stochastic}, Theorem~8.13. Using \eqref
{eq:lyap_CTMC} and the nonnegativity of $V$, we have for all
sufficiently large $n$ and all $t>0$ that
%
\begin{equation}
\frac{1}{t}\Ex_x \biggl[ \int_0^t
V\bigl(\widehat{X}^n(s)\bigr) \,ds \biggr]\leq\frac{1}{t}2\delta
^{-1}\bigl(V(x)+b t \bigr).
\end{equation}
Letting $\nu^n$ be the steady-state distribution of $\hX^n$ we have,
for each $M$, that
\[
\Ex_{\nu^n}\bigl[V\bigl(\widehat{X}^n(0)\bigr)\wedge M\bigr]=
\lim_{t\rightarrow\infty
}\frac{1}{t}\Ex _x \biggl[\int
_0^t V\bigl(\widehat{X}^n(s)\bigr)
\wedge M \,ds \biggr]\leq2\delta^{-1}b.
\]
The result now follows from the nonnegativity of $V$ and the monotone
convergence theorem.
\end{pf*}

\section{Two examples}\label{sec:examples}

Lyapunov functions that satisfy Assumption~\ref{asum:L} must be
identified on a case-by-case basis. For the first example---the
Erlang-A queue---this is a straightforward task. For the second
example---a queue with many servers and phase-type service time
distribution---this task is substantially more difficult, but recent
work \cite{dieker2011positive} provides us with the required function.

\subsection{The Erlang-A queue}

We consider a sequence of queues with a single pool of i.i.d. servers
that serve one class of impatient i.i.d. customers. Arrivals follow a
Poisson process (with rate $n$ in the $n$th queue), service times are
exponentially distributed with rate $\mu$ and customers' patience times
are exponentially distributed with rate $\theta$. In the $n$th queue,
there are $N^n$ servers in the server pool. Let $X^n(t)$ be the total
number of jobs in the $n$th queue (waiting or in service) at time $t$.
Then $(X^n(t),t\geq0)$ is a birth and death process with state space
$\bbZ_+$, birth rate $n$ in all states and death rate $\mu(x\wedge
N^n)+\theta(x-N^n)^+$ in state $x$ where, for the remainder of the
paper, we use $(x)^+=\max\{0,x\}$, $(x)^{-}=\max\{0,-x\}$.
\textit{We assume that $\theta>0$} so that positive recurrence of $X^n$
follows easily.

The drift $F^n$ is then specified here by
\[
F^n(x)=n-\mu\bigl(x\wedge N^n\bigr)-\theta
\bigl(x-N^n\bigr)^+,\qquad x\in\bbZ_+,
\]
and is trivially extended here to the real line by allowing $x$ to take
real values (including negative values). The FM is then given by
\renewcommand{\theequation}{FM}
\begin{equation}
 \bar{x}^n(t)=\bar{x}^n(0)+\int
_0^t F^n\bigl(\bar
{x}^n(s)\bigr)\,ds.
\end{equation}
There exists a unique point $\bar{x}^n_{\infty}$ in which $F^n(\bar
{x}^n_{\infty})=0$. At this point
$n=\mu(\bar{x}^n_{\infty}\wedge N^n)+\theta(\bar{x}^n_{\infty}-N^n)^+$
so that
\[
\bar{a}^n(0)= \frac{1}{n} \bigl(n+\mu\bigl(
\bar{x}^n_{\infty}\wedge N^n\bigr)+\theta\bigl(\bar
{x}^n_{\infty}-N^n\bigr)^+ \bigr)\equiv2.
\]
The DM for the Erlang-A queue is subsequently given by
\renewcommand{\theequation}{DM}
\begin{equation}
\hY^n(t)=\hY^n(0)+\int_0^t
\hF^n\bigl(\hY^n(s)\bigr)\,ds +\sqrt{2}B(t),\vadjust{\goodbreak}
\end{equation}
where
\[
\hF^n(x)=\mu \bigl(\bigl(f^n(x)\bigr)^{-}-
\bigl(f^n(0)\bigr)^{-} \bigr)-\theta \bigl(
\bigl(f^n(x)\bigr)^+-\bigl(f^n(0)\bigr)^+ \bigr),
\]
and $f^n(x)=x+(\bar{x}^n_{\infty}-N^n)/\sqrt{n}$. It is easily verified
that there exists $\eta>0$ such that $\hF^n(x)\leq-\eta x$ when $x>0$
and $\hF^n(x)\geq-\eta x $ if $x<0$. Fixing $\varrho\geq1$ and taking
\[
V_m(x)=\varrho+x^{2m},\qquad x\in\bbR, m\in\bbN,
\]
we have that $V_m(x)>V_m(0)$ for all $x\neq0$ and
\[
DV_m(x)\hF^n(x)\leq-\eta(2m) \bigl(V_m(x)-V_m(0)
\bigr)\qquad\mbox{for all } x\neq0.
\]
Note that
$V_m$ is trivially sub-exponential.
Further, for all sufficiently large $|x|$,
\[
D^2 V_m(x)= 2m(2m-1)x^{2m-2}\leq
\frac{\eta}{2} x^{2m},
\]
so that the conditions of Lemma~\ref{lem:FMtoDM} are satisfied and, in
turn, UL holds for the~DM. Further, for each $t\geq0$,
$X^n(t)\leq X^n(0)+ N^n+A^n(t)$ where $A^n(t)$ is the number of
arrivals by time $t$. Condition \eqref{eq:finite_integral} then follows
from basic properties of the Poisson process. We have the following consequence.

\begin{lem} Fix $\varrho\geq1$ and positive $m\in\bbN$. Then,
$V_m(x)=\varrho+x^{2m}$ satisfies Assumption~\ref{asum:L} for the DM of
the Erlang-A queue.
\end{lem}

Fixing $m\in\bbN$ and choosing sufficiently large $\varrho$, we can
take $\bar{V}_m=V_{4m}$ in Corollary~\ref{cor:tightness}; see Remark~\ref{rem:simplecase}. The following is now a direct consequence of
Theorem~\ref{thmm:main} and Corollary~\ref{cor:tightness}.

\begin{thmm}[(Approximation gap for the Erlang-A queue in stationarity)]\label{thmm:ErlangA}
Consider a sequence of Erlang-A queues as above and let $f$ be such
that $\bar{f}\leq V_m$ for some $m\in\bbN$. Then
\[
\limsup_{n\rightarrow\infty } \nu^n\bigl(|f|\bigr)<\infty\quad\mbox{and} \quad\nu
^n(f)-\pi ^n(f)=\calO (1/\sqrt{n}).
\]
\end{thmm}

\begin{rem}[(Universality and the connection to \cite{Excursion})]
Above, we did not impose any restrictions on the way in which the
number of servers, $N^n$, scales with $n$ so that one may interpret our
DM as a \textit{universal} approximation for the Erlang-A queue.
Universality for this queue (and its contrast with the assumption of a
so-called operational regime) are discussed at length in \cite
{Excursion}; see also the references therein. A similar result is
proved there for the Erlang-A queue using an approach that, while
having important similarities to the approach we take here, is based on
approximating the excursions of the process $X^n$ above and below $\bar
{x}_{\infty}^n$. In this one-dimensional Markov chain, the Poisson
equation we use here is (informally) a ``pasting'' of the Dirichlet
problems studied in \cite{Excursion}.

In their greatest generality, the results of \cite{Excursion} are not a
special case of Theorem~\ref{thmm:ErlangA} above. In \cite{Excursion}
the authors allow the service rate $\mu$ to vary with $n$. This is
facilitated by the excursion approach taken there but violates the
assumptions required to apply our results, particulary, the uniform
Lipschitz continuity of $\hF^n$. Moreover, the approach in \cite
{Excursion} seems to be easily extendable to the case with $\theta=0$
in which case the DM is not exponentially ergodic and Assumption~\ref
{asum:L} is not satisfied.
\end{rem}

\subsection{A phase-type queue with many servers}
\label{sec:PH}
We next consider the single class $M/PH/n+M$ queue. This is a
generalization of the Erlang-A queue where the exponential service time
is replaced by a phase-type service-time; see \cite{dai2010many} for a
detailed construction. We repeat here only the essential details.

Let $I$ be the number of service phases, and let $1/\nu_k$ be the
average length of phase $k=1,\ldots,I$. We assume that $p=(1,\ldots
,0)'$, corresponding to all customers commencing their service at phase
1 (the diffusion limits in \cite{dai2010many} cover the general case
where $p$ is an arbitrary probability vector). Having completed phase~$i$
 a job transitions into phase $j$ with probability $P_{ij}$. The
triplet $(p,\nu,P)$ defines the phase-type service-time distribution.

Let
\[
R=\bigl(I-P'\bigr)\operatorname{diag}(\nu)\quad\mbox{and}\quad 1/\mu=e'R^{-1}p,\qquad
\gamma=\mu R^{-1}p.
\]
Note that $\sum_{k}\gamma_k=1$. As before, the patience rate is
$\theta>0$.

We consider a sequence of such queues indexed by the \textit{arrival rate
$n\in\bbZ_+$}. Let
\[
\gamma^n=n\gamma,\qquad n\in\bbN.
\]
Let $X_1^n(t)$
be the number of customers in the first phase of their service and
waiting in the queue at time $t$. For $i>1$, let $X_i^n(t)$ be the
number of customers in phase $i$ of service at time $t$. The process
\[
X^n(t)=\bigl(X_1^n(t),\ldots,X_I^n(t)
\bigr),
\]
is then a CTMC.

For simplicity of exposition we assume here that $\sum_{k}\gamma_k^n$
is integer valued for each $n$ and that the number of servers $N^n$
satisfies $N^n=\sum_{k}\gamma_k^n$. This implies, trivially, that
$N^n=\sum_{k}\gamma_k^n+\calO(\sqrt{n})$ which corresponds to the
so-called Halfin--Whitt many-server regime and allows us subsequently
to build on the results of~\cite{dai2010many} and \cite
{dieker2011positive} that study diffusion limits in this regime. The
analysis below is easily extended to the case $N^n=\sum_{k}\gamma
_k^n+\beta\sqrt{n}+o(\sqrt{n})$ for some $\beta\neq0$.

Define
\[
\bar{x}^n_{\infty}=\bigl(\gamma_1^n,
\ldots,\gamma_I^n\bigr),
\]
and the scaled and centered process $\hX^n$ as in \eqref{eq:hxdefin}. Then,
%
\renewcommand{\theequation}{\arabic{equation}}
\setcounter{equation}{47}
\begin{equation}\label{eq:PH_drift}
\hF_i^n(x)=\cases{ %
\displaystyle-
\nu_i x_i+\sum_{k\neq i,k\neq1}P_{ki}
\nu_k x_k+ \nu _1P_{1i}
\bigl(x_1-\bigl(e'x\bigr)^+\bigr),&\quad$\mbox{if } i\neq1,$
\vspace*{2pt}\cr
\displaystyle -\nu_1 \bigl(x_1-\bigl(e'x\bigr)^+
\bigr)-\theta\bigl(e'x\bigr)^+,&\quad $\mbox{if } i=1.$}\hspace*{-30pt}
\end{equation}
This is written, in Matrix notation, as
%
\begin{eqnarray}
\qquad\hF^n(x)&=&-R x+(R-\theta I)p\bigl(e'x
\bigr)^+.\label{eq:PH_drift_matrix}
\\
\label{eq:PH_diff1}
n\bar{a}_{kk}^n(x)&=&\cases{ %
\displaystyle\sum_{i\neq k,i\neq1}P_{ik}\nu_k \bigl(
\gamma _k^n+\sqrt {n} x_k\bigr)+
\nu_k \bigl(\gamma_k^n+\sqrt{n}
x_k\bigr) &
\vspace*{2pt}\cr
\quad{}\displaystyle\hspace*{18pt} +\nu_1P_{1k} \bigl(\gamma _1^n+
\sqrt{n} x_1-\sqrt{n}\bigl(e'x\bigr)^+\bigr),&\quad $\mbox{if
} k\neq1,$
\vspace*{2pt}\cr
n+\nu_1\bigl(\gamma_1^n+\sqrt{n}
x_1\bigr)+\theta\sqrt{n}\bigl(e'x\bigr)^+,&\quad $\mbox{if
} k=1,$}
\end{eqnarray}
and, for $k\neq j$,
%
\begin{equation}\label{eq:PH_diff2}
\quad n\bar{a}_{kj}^n(x)=\cases{ %
P_{kj}\nu_k \bigl(\gamma_k^n+
\sqrt{n} x_k\bigr)+P_{jk}\nu_j \bigl(
\gamma_j^n+\sqrt{n} x_j\bigr),&\quad $\mbox{if } k
\neq1,$
\vspace*{2pt}\cr
P_{kj}\nu_k \bigl(\gamma_k^n+
\sqrt{n}x_k-\sqrt{n}\bigl(e'x\bigr)^+
\bigr)\vspace*{2pt}\cr
\qquad{}+P_{jk}\nu_j \bigl(\gamma _j^n+
\sqrt{n} x_j\bigr),&\quad $\mbox{ if } k=1.$}
\end{equation}

The functions $\hF^n$ and $\bar{a}^n$ satisfy \eqref{eq:Flip} and
\eqref{eq:alip}. Assumption~\ref{asum:base} holds in this example as
the chain is trivially nonexplosive and irreducible. The positive
recurrence follows immediately from the fact that $\theta>0$.

The diffusion model is given by
\renewcommand{\theequation}{DM}
\begin{equation}
\hY^n(t)=y+\int_0^t
\hF^n\bigl(\hY^n(s)\bigr)\,ds+\sqrt{\bar{a}^n(0)}
B(t),
\end{equation}
with $\hF^n$ as in \eqref{eq:PH_drift_matrix} and diffusion coefficient
$\bar{a}^n$ as in \eqref{eq:PH_diff1}--\eqref{eq:PH_diff2}. Note
\eqref
{eq:PH_drift_matrix}--\eqref{eq:PH_diff2} that $\hF^n$ and $\bar
{a}^n(0)$ do not, in fact, depend here on $n$. The existence of a
quadratic Lyapunov function, $V$, for $\hY^n$ then follows from \cite{dieker2011positive},
Theorem~3---this function is specified in
equation (5.24) there. (To extend this argument to the general case
with $N^n=\sum_{k}\gamma_k^n+\beta\sqrt{n}+o(\sqrt{n})$, note that $V$
in \cite{dieker2011positive} is still a Lyapunov function for each $n$
if we perturb $\hF^n$ by a constant and $\bar{a}^n(0)$ by a term that
shrinks proportional to $1/\sqrt{n}$.)

With a careful choice of the smoothing function $\phi$ there, the
function $\Psi=\varrho+V$ (for any constant $\varrho\geq1$) is also
sub-exponential. Finally, \eqref{eq:finite_integral} is argued as in
the Erlang-A case using crude bounds on the Poisson arrivals.

The function $\Psi=\varrho+V$ thus satisfies Assumption~\ref{asum:L}.
It is easily verified that $\Psi\in\mC^3(\bbR^d$) and satisfies
\eqref
{eq:3Diff1} so that, as in Remark~\ref{rem:simplecase}, $\Psi
_m(x)=(\Psi
(x))^m$ satisfies Assumption~\ref{asum:L} with re-defined constants
$\delta,b$ and $K$. Choosing sufficiently large $\varrho$ guarantees
that $\Psi_{4m}(\cdot)\geq\Psi_m(\cdot)(1+|\cdot|)^4$. The following
is then an immediate consequence of Theorem~\ref{thmm:main} and
Corollary~\ref{cor:tightness}.

\begin{cor} Consider the sequence of phase-type queues as above, and
let $f$ be such that $\bar{f}\leq\Psi_m$ for some $m\in\bbN$. Then
\[
\limsup_{n\rightarrow\infty } \nu^n\bigl(|f|\bigr) <\infty\quad\mbox{and}\quad\nu
^n(f)-\pi ^n(f) =\calO (1/\sqrt{n}).
\]
\end{cor}

Thus, as in Remark~\ref{rem:AllInOne}, we have a Lyapunov function that
allows us to establish simultaneously the stability of the Markov chain
for each sufficiently large $n$, the uniform integrability of moments
and the approximation gap. It is worth noting that the fact that
$\limsup_{n\rightarrow\infty } \nu^n(|f|) <\infty$ was already
established, by
alternative means and for more general (multiclass) phase-type queues,
in \cite{DaiDiekerGao}.

\section{Proof of Theorem \texorpdfstring{\protect\ref{thmm:expo_ergo}}{3.3}}
\label{sec:expo_ergo}

The main step in this proof is a uniform minorization condition for a
time-discretized version of $\hY^n$. Once this is established (see
Lemma~\ref{lem:minorization} below), we build on \cite
{baxendale2005renewal} to complete the argument. The proofs of the
lemmas that are stated in this section appear in the \hyperref[app]{Appendix}.

We first consider a linear transformation of $\hY^n$. Specifically, let
$L_n$ be the unique square root of the matrix $\bar{a}^n(0)$; see
\cite{HoJ:94}, Theorem~7.2.6. In particular, $L_n(L_n)^T =\bar{a}^n(0)$. The
matrix $L_n$ is itself invertible and its inverse is the square root of
the inverse of $\bar{a}^n(0)$; see \cite{HoJ:94}, page 406. Let
%
\renewcommand{\theequation}{\arabic{equation}}
\setcounter{equation}{51}
\begin{equation}
\hF _L^n(x)=L_n^{-1}
\hF^n(L_n x), \qquad x\in\bbR^d,\label{eq:hFL}
\end{equation}
and define
\[
Z_L^n(t)=L_n^{-1}
\hY^n(t),\qquad t\geq0.
\]
Then $Z_L^n$ is a $d$-dimensional Brownian motion with drift $\hF^n$,
that is,
\[
Z_L^n(t)=z+\int_0^t
\hF_L^n\bigl(Z_L^n(s)
\bigr)\,ds+B(t),
\]
where $z =L_n^{-1}\hY^n(0)$.

We next consider the discrete-time analogues of both $Z_L^n$ and $\hY
^n$. Let
\[
\Phi^n_l=Z_L^n(l)\quad \mbox{and} \quad\psi^n_l=\hY^n(l)\qquad\mbox{for } l\in
\bbZ_+.
\]
Let $\Pd_{\Phi^n}(\cdot,\cdot)$ and $\Pd_{\psi^n}(\cdot,\cdot)$
be the
corresponding one-step transition functions. Below $\mathcal{B}(\bbR
^d)$ is the family of Borel sets in $\bbR^d$.

\begin{lem}\label{lem:minorization} Fixing $K>0$, there exist a probability measure $\mathcal
{Q}$ with $\mathcal{Q}(B_0(K))=1$ and a constant $\epsilon<1$ (both not
depending on $n$) such that
\[
\Pd_{\Phi^n}(x,\mathcal{E})\geq\epsilon\mathcal{Q}(L_n
\mathcal {E} ),\qquad x\in L_n^{-1} B_0(K),
\mathcal{E}\in\mathcal{B}\bigl(\bbR^d\bigr).
\]
There consequently exists a constant $\widetilde{\epsilon}<1$ (not
depending on $n$) such that
\[
\Pd_{\psi^n}(x,\mathcal{E})\geq\widetilde{\epsilon}\mathcal {Q}(\mathcal
{E}),\qquad x\in B_0(K), \mathcal{E}\in\mathcal{B}\bigl(\bbR^d
\bigr).
\]
\end{lem}

The following translates the Lyapunov property UL into the discrete
time setting.

\begin{lem}\label{lem:fromdisctocont} Let $V$ be as in Assumption~\ref{asum:L}. Then there exist
finite positive constants $\gamma<1$ and $\bar{b}$ (not depending on
$n$) such that for all $n\in\bbN$ and all $x\in\mathbb{R}^d$,
%
\begin{equation}
\Ex _x\bigl[V\bigl(\hY^n(1)\bigr)\bigr]\leq(1-
\gamma)V(x)+\bar{b}\mathbh{1}_{\overline
{B}_0(K)}(x).\label{eq:discretedrift}
\end{equation}
\end{lem}

Using the fact that $V(x)\rightarrow\infty $ as $|x|\rightarrow
\infty $, \eqref{eq:discretedrift}
implies that there exist finite positive constants $K$, $\lambda<1$ and
$M$ such that
%
\begin{equation}
\label{eq:baxendalcond} \Ex_x\bigl[V\bigl(\hY^n(1)\bigr)\bigr]
\leq \cases{ %
\lambda V(x),&\quad $\mbox{if } x\notin
\overline {B}_0(K),$
\vspace*{2pt}\cr
M,& \quad$\mbox{if } x\in\overline{B}_0(K).$}
\end{equation}

The following is then a direct consequence of \cite{baxendale2005renewal},
Theorem~1.1. Assumptions~(A1)--(A3)  there hold by
Lemmas \ref{lem:minorization}, \ref{lem:fromdisctocont} and by \eqref
{eq:baxendalcond}.

\begin{cor} There exist constants $\mathbb{M}$ and $\mu$ (not depending
on $n$) such that for each $m\in\bbN$,
\[
\sup_{n}\sup_{x\in\bbR^d}\sup
_{|f|\leq V}\frac{1}{V(x)}\bigl |\Ex _x\bigl[f(\hY
^n(m)\bigr]-\pi^n(f) \bigr|\leq\mathbb{M}e^{-\mu m}.
\]
\end{cor}

With these we are ready for the proof of Theorem~\ref{thmm:expo_ergo}.
\begin{pf*}{Proof of Theorem \protect\ref{thmm:expo_ergo}} The proof
of the
theorem now follows as in \cite{TM:93}, page 536.
Specifically, let $s=t-\lfloor t\rfloor$
\begin{eqnarray*}
\sup_{|f|\leq V}\bigl|\Ex_x\bigl[f\bigl(\hY^n(t)
\bigr)\bigr]-\pi^n(f)\bigr|& =& \sup_{|f|\leq V}\bigl|
\Ex_x\bigl[f\bigl(\hY^n\bigl(\lfloor t\rfloor+s\bigr)
\bigr)\bigr]-\pi^n(f)\bigr|
\\
& =&\sup_{|f|\leq V}\bigl|\Pd_{\hY^n}^s(x,dy) \bigl(
\Ex_x\bigl[f\bigl(\hY^n\bigl(\lfloor t\rfloor \bigr)
\bigr)\bigr]-\pi^n(f)\bigr)\bigr|
\\
&\leq&\int_{y}\Pd_{\hY^n}^s(x,dy)\sup
_{|f|\leq
V}\bigl|\Ex _y\bigl[f\bigl(\hY^n\bigl(
\lfloor t\rfloor\bigr)\bigr)\bigr]-\pi^n (f)\bigr|
\\
& \leq& \mathbb{M}e^{-\mu\lfloor t\rfloor} \Ex_x\bigl[V\bigl(
\hY^n(s)\bigr)\bigr]
\\
&\leq& \mathbb{M} e^{\mu} e^{-\mu t}
\bigl(V(x)+b\bigr),
\end{eqnarray*}
where $\Pd_{\hY^n}^s(x,A)$ is the transition probability function of
$\hY^n$ in $s$ time units. In the last inequality we used \eqref
{eq:tbound} and the fact that $s=t-\lfloor t\rfloor\leq1$. Finally,
since $V\geq1$, the theorem holds with the constants $\mathcal
{M}=\mathbb{M}e^{+1}(1+b)$ and $\mu$.
\end{pf*}
\section{Concluding remarks}
\label{sec:conclusions} Diffusion models
are useful in the approximation of Markov chains. We proved that, under
a uniform Lyapunov condition, the steady-state of some multidimensional
CTMCs can be approximated with impressive accuracy by the steady-state
of a relatively tractable diffusion \textit{model}.

The existence of a diffusion limit that satisfies the Lyapunov
requirement---as is the case for the phase-type queue considered in
Section~\ref{sec:PH}---can facilitate the application of our results.
The distinction between the diffusion model and diffusion limit is,
however, important. A central motivation behind this work is to bypass
the need for diffusion limits with the objective of providing
steady-state diffusion\vadjust{\goodbreak} approximation whose precision does not depend on
assumption with regards to limiting values of underlying parameters.
That is, we ultimately seek to provide ``limit-free'' (or \emph{universal}) approximations.

A uniform Lyapunov condition, as we require in Assumption~\ref{asum:L},
need not hold in general. Informally, one expects such a condition to
hold if the scale parameter $n$ has limited effect on the drift of the
process around the FMs stationary point. Many-server queues with
abandonment, as those we use to illustrate our results, seem to satisfy
this characterizations: diffusion limits (regardless of the parameter
regime, determining how the number of servers $N^n$ scales with $n$)
are generalizations of the OU process. It remains to identify the
broadest characterization of Markov chains for which a uniform Lyapunov
condition can be expected to hold.

In addition, the following extensions seem important:

\textit{State-space collapse}. A fundamental phenomenon in
diffusion limits for multi-class queueing systems is that of
state-space collapse (SSC). With SSC, the diffusion limit ``lives'' on
a state-space that is of lower dimension relative to the \mbox{original} CTMC:
some coordinates of the CTMC become, asymptotically, deterministic
functions of others. For example, if one allows for arbitrary
initial-phase vectors $p$ in the example of Section~\ref{sec:PH}, the
number of customers in queue with initial phase $k$ is asymptotically
equal to $p_k$; see \cite{dai2010many}. To exploit state-space collapse
within the diffusion-model framework used in this paper, one must
develop bounds (rather than convergence results) for state-space collapse.

\textit{Single server queues and reflection.} A key challenge
with single-server queueing systems is that of reflection. Such
reflection may violate our assumptions on $\hF^n$. Consider, for
example, the $M/M/1+M$ queue---this is a single-server version of the
Erlang-A queue discussed in Section~\ref{sec:examples}. Suppose that
the arrival rate and service rate in the $n$th queue satisfy $\lambda
^n=n\lambda$, and $\mu^n=\lambda^n-\beta\sqrt{n}$ (for $\beta
>0$). Let
$\theta>0$ be the patience parameter. Then
\begin{eqnarray*}
F^n(x)&=&\lambda^n -\mu^n \mathbh{1}\{x>0\}-
\theta x
\\
& =&\beta\sqrt{n}-\theta x+\mu^n\mathbh{1}\{x=0\},
\end{eqnarray*}
so that $\bar{x}^n_{\infty}=\beta\sqrt{n}/\theta$. Also, $\hF
^n(-\beta
/\theta)=F^n(0)/\sqrt{n}=\beta+\mu^n/\sqrt{n}$ and, in particular
$|\hF
^n(-\beta/\theta)-\hF^n(0)|=\beta+\mu^n/\sqrt{n}=\sqrt{n}\lambda
\rightarrow\infty$ as $n\rightarrow\infty $. Clearly, \eqref
{eq:Flip} is violated.

It is fair to conjecture that similar results as ours can be proved in
such settings provided that the reflection is explicitly captured in
the DM. Extending our results to DMs with reflection seems to present a
challenge insofar as the theory of PDEs that arise from the Poisson
equation for such networks is less developed and poses a challenge in
terms of the gradient bounds that are central to our analysis here;
see, for example, \cite{budhiraja2007long}, where the Poisson equation
for constrained diffusion is discussed as well as, in the context of
ergodic control, \cite{borkar2004ergodic}.

\begin{appendix}\label{app}
\section*{Appendix}
\begin{pf*}{Proof of Lemma \protect\ref{lem:ito}} Fix $x\in\widehat{E}^n$.
By It\^{o}'s rule applied to the pure jump process $(\Psi(\widehat
{X}^n(t)), t\geq
0)$ we have that
%
\begin{eqnarray}
\Psi\bigl(\widehat{X}^n(t)\bigr)&=&\Psi(x)+\sum
_{s\leq t}\sum_{i=1}^d \Psi
_i\bigl(\widehat{X}^n (s-)\bigr)\Delta
\widehat{X}^n_i(s)
\nonumber
\\[-8pt]
\\[-8pt]
\nonumber
&&{}+ \sum_{s\leq t} \Biggl[\Psi\bigl(
\widehat{X}^n(s)\bigr)-\Psi \bigl(\widehat{X}^n (s-)
\bigr)-\sum_{i=1}^d \Psi_i
\bigl(\widehat{X}^n(s-)\bigr)\Delta\widehat{X}^n_i(s)
\Biggr]. \label{eq:itojump}
\end{eqnarray}
From the linear growth of $\hF^n$ and from \eqref
{eq:integrability_cond}, it then follows that
\[
\Ex_x \biggl[\int_0^t \bigl|D\Psi
\bigl(\widehat{X}^n(s)\bigr)\bigr|\bigl|\hF^n\bigl(\widehat
{X}^n(s)\bigr)\bigr|\,ds \biggr]<\infty.
\]
We can then apply L\'evy's formula for CTMCs (see, e.g., \cite{Bre:80}, Exercise
I.2.E2) to get that
\[
\sum_{s\leq t}\sum_{i=1}^d
\Psi_i\bigl(\widehat {X}^n(s-)\bigr)\Delta
\widehat{X}^n _i(s)-\sum_{i=1}^d
\int_0^t \Psi_i\bigl(
\widehat{X}^n(s)\bigr)\hF _i^n\bigl(
\widehat{X}^n(s)\bigr)\,ds
\]
is a martingale with respect to the filtration in \eqref{eq:filtration}
and, in turn, for all $t\geq0$,
\[
\Ex_x \Biggl[\sum_{s\leq t}\sum
_{i=1}^d \Psi_i\bigl(\widehat
{X}^n(s-)\bigr)\Delta\widehat{X}^n _i(s)
\Biggr]=\Ex_x \Biggl[\sum_{i=1}^d
\int_0^t \Psi_i\bigl(\widehat
{X}^n(s)\bigr)\hF _i^n\bigl(
\widehat{X}^n(s)\bigr)\,ds \Biggr].
\]
To treat the second line of \eqref{eq:itojump}, we decompose it into
%
\renewcommand{\theequation}{$\mathrm{D}$}
\begin{equation}\label{eqD}
\frac{1}{2}\sum_{s\leq t}\sum
_{i,j}^d \Psi_{ij}\bigl(\widehat
{X}^n(s-)\bigr)\Delta \widehat{X}^n _i(s)
\Delta\widehat{X}^n_j(s)
\end{equation}
and
%
\renewcommand{\theequation}{$\mathrm{A}$}
\begin{eqnarray}\label{eqA}
&&
\sum_{s\leq t} \Biggl[\Psi\bigl(
\widehat{X}^n(s)\bigr)-\Psi\bigl(\widehat {X}^n(s-)\bigr)-
\sum_{i=1}^d \Psi _i\bigl(
\widehat{X}^n(s-)\bigr)\Delta\widehat{X}^n_i(s)
\nonumber
\\[-8pt]
\\[-8pt]
\nonumber
&&\hspace*{60pt}\qquad{} -\frac{1}{2}\sum_{i,j}^d
\Psi_{ij}\bigl(\widehat{X}^n (s-)\bigr)\Delta
\widehat{X}^n_i(s)\Delta\widehat{X}^n_j(s)
\Biggr].
\end{eqnarray}
We treat (\ref{eqD}) first. By \eqref{eq:alip}, $|\bar{a}^n(x)|\leq
|\bar
{a}^n(0)|+K_a|x|/\sqrt{n}$ so that, by \eqref{eq:integrability_cond},
\[
\Ex_x \biggl[\int_0^t
\bigl|D^2\Psi\bigl(\widehat{X}^n(s)\bigr)\bigr|\bigl |\bar
{a}^n\bigl(\widehat{X}^n(s)\bigr)\bigr|\,ds \biggr] <\infty ,\qquad t
\geq0, x\in\widehat{E}^n,
\]
and applying L\'evy's formula once again, we obtain
\begin{eqnarray*}
&&\frac{1}{2}\Ex_x \Biggl[\sum_{s\leq t}
\sum_{i,j}^d \Psi _{ij}\bigl(
\widehat{X}^n (s-)\bigr)\Delta\widehat{X}^n_i(s)
\Delta\widehat{X}^n_j(s) \Biggr]
\\
&&\qquad = \frac{1}{2}\Ex_x \Biggl[\sum
_{i,j}^d\sum_{\ell
}\int
_0^t \Psi_{ij}\bigl(
\widehat{X}^n(s)\bigr)\ell_i\ell_j
\frac{1}{n}\beta _{\ell
}^n\bigl(\sqrt {n}
\widehat{X}^n(s)+\bar{x}^n_{\infty}\bigr)\,ds \Biggr]
\\
&&\qquad =\frac{1}{2} \Ex_x \Biggl[\sum
_{i,j}^d \int_0^t
\Psi_{ij}\bigl(\widehat{X}^n(s)\bigr)\bar{a}_{ij}^n
\bigl(\widehat {X}^n(s)\bigr)\,ds \Biggr]
\\
&&\qquad = \frac{1}{2}\Ex_x \Biggl[\sum
_{i,j}^d \int_0^t
\Psi_{ij}\bigl(\widehat{X}^n (s)\bigr)\bar
{a}_{ij}^n(0)\,ds \Biggr]
\\
&&\qquad\quad{} +\frac{1}{2}\Ex_x \Biggl[\sum
_{i,j}^d \int_0^t
\Psi_{ij}\bigl(\widehat{X}^n(s)\bigr) \bigl(
\bar{a}^n\bigl(\widehat {X}^n(s)\bigr)-\bar
{a}_{ij}^n(0)\bigr)\,ds \Biggr].
\end{eqnarray*}
The second item in the last line is $D_{\Psi}^{n,x}(t)$ in the
statement of the lemma. We have proven thus far that
\begin{eqnarray*}
&& \Ex_x\bigl[\Psi\bigl(\widehat{X}^n(t)\bigr)
\bigr]
\\
&&\qquad = \Psi(x)+\Ex_x \Biggl[\sum_{i=1}^d
\int_0^t \Psi_i\bigl(\widehat
{X}^n(s)\bigr)\hF _i^n\bigl(
\widehat{X}^n (s)\bigr)\,ds \Biggr]
\\
&&\quad\qquad{}+ \frac{1}{2}
\Ex_x \Biggl[\sum_{i,j}^d \int
_0^t \Psi_{ij}\bigl(
\widehat{X}^n (s)\bigr)\bar {a}_{ij}^n(0)\,ds
\Biggr]
+D_{\Psi}^{n,x}(t)+A_{\Psi
}^{n,x}(t)
\\
&&\qquad  = \Psi(x)+\Ex_x \biggl[\int_0^t
\mathcal{A}^n\Psi\bigl(\widehat{X}^n (s)\bigr)\,ds
\biggr]+D_{\Psi}^{n,x}(t)+A_{\Psi}^{n,x}(t),
\end{eqnarray*}
where $D_{\Psi}^{n,x}$ is as in the statement of the lemma and
$A_{\Psi
}^{n,x}(t)=\Ex_x[\mathbf{A}]$ (we will prove below that this expectation
is well defined). To bound $A_{\Psi}^{n,x}$ note that, by Taylor's theorem,
\begin{eqnarray*}
&&\Psi\bigl(\widehat{X}^n(s)\bigr)-\Psi\bigl(\widehat{X}^n(s-)
\bigr)\\
&&\qquad=\sum_{i=1}^d \Psi _i
\bigl(\widehat{X}^n (s-)\bigr)\Delta\widehat{X}^n
_i(s)
\\
&&\qquad\quad{}+\frac{1}{2}\sum_{i,j}^d
\Psi_{ij}\bigl(\widehat{X}^n(s-)+\eta _{\widehat{X}^n
(s-),\widehat{X}^n
(s)}\bigr)
\Delta\widehat{X}^n_i(s)\Delta\widehat{X}^n_j(s),
\end{eqnarray*}
where $\eta_{\widehat{X}^n(s-),\widehat{X}^n(s)}\in\prod
_{i=1}^d[\widehat{X}^n_i(s-),\widehat{X}^n
_i(s)]$. Thus
\[
\mathbf{A}=\frac{1}{2} \sum_{s\leq t} \Biggl( \sum
_{i,j}^d \bigl(\Psi _{ij}\bigl(
\widehat{X}^n(s-)+\eta_{\widehat{X}^n(s-)}\bigr)-\Psi _{ij}
\bigl(\widehat{X}^n(s-)\bigr) \bigr)\Delta \widehat{X}^n
_i(s)\Delta\widehat{X}^n_j(s) \Biggr).
\]
Here note that $|\Delta\hX_i^n(s)||\Delta\hX_j^n(s)|\leq\bar{\ell
}^2/n$. Let $\widetilde{\Psi}_{ij}(x,y)=\Psi_{ij}(x+\eta
_{x,y})-\Psi
_{ij}(x)$. Note that $|\widetilde{\Psi}_{ij}(x,y)|\leq\frac{\bar
{\ell
}}{\sqrt{n}}[\Psi]_{2,1,B_{x}(\bar{\ell}/\sqrt{n})}$ for $x,y\in
\hE^n$
with $y\in\break  B_x(\bar{\ell}/\sqrt{n})$. Since $\sum_{\ell}|\ell
_i||\ell
_j|\beta_{\ell}^n(x)\leq\sum_{\ell}(|\ell_i|^2+|\ell_j|^2)\beta
_{\ell
}^n(x)\leq|a^n(x)|$, we have that
\begin{eqnarray*}
&& \frac{1}{2 n}\Ex_x \Biggl[\sum
_{i,j}^d \int_0^t
\sum_{\ell} \bigl|\widetilde {\Psi}_{ij}\bigl(
\widehat{X}^n(s),\widehat{X}^n(s)+\ell/\sqrt{n}\bigr)\bigr| |
\ell_i||\ell _j|\beta_{\ell
}^n
\bigl(X^n(s)\bigr)\,ds \Biggr]
\\
&& \qquad\leq \frac{\bar{\ell}}{\sqrt{n}}\frac{1}{2 n}\Ex_x \Biggl[\sum
_{i,j}^d \int_0^t
[\Psi]_{2,1,B_{\widehat{X}^n(s)}(\bar{\ell}/\sqrt{n})}\sum_{\ell
} |\ell
_i||\ell _j|\beta_{\ell}^n
\bigl(X^n(s)\bigr)\,ds \Biggr]
\\
&&\qquad \leq \frac{\bar{\ell}}{\sqrt{n}}\frac{1}{2 n}\Ex_x \biggl[\int
_0^t [\Psi ]_{2,1,B_{\widehat{X}^n(s)}(\bar{\ell}/\sqrt
{n})}\bigl|a^n
\bigl(X^n(s)\bigr)\bigr|\,ds \biggr]
\\
&&\qquad = \frac{\bar{\ell}}{2\sqrt{n}}\Ex_x \biggl[\int_0^t
[\Psi ]_{2,1,B_{\widehat{X}^n(s)}(\bar{\ell}/\sqrt{n})}\bigl|\bar {a}^n\bigl(\widehat{X}^n (s)
\bigr)\bigr|\,ds \biggr] <\infty,
\end{eqnarray*}
where the finiteness follows from \eqref{eq:alip} and condition \eqref
{eq:integrability_cond}.

We can apply L\'evy's formula one final time to conclude that
\begin{eqnarray*}
\bigl|\Ex_x[\mathbf{A}]\bigr|&=& \Biggl\llvert \frac{1}{2 n}\sum
_{i,j}^d \Ex_x \biggl[ \int
_0^t \sum_{\ell}
\widetilde{\Psi}_{ij}\bigl(\widehat{X}^n(s),\widehat
{X}^n(s)+\ell/\sqrt {n}\bigr)\ell _i\ell_j
\beta_{\ell}^n\bigl(X^n(s)\bigr)\,ds\biggr]\Biggr
\rrvert
\\
& \leq&\frac{\bar{\ell}}{2\sqrt{n}}\Ex_x \biggl[\int_0^t
[\Psi ]_{2,1,B_{\widehat{X}^n(s)}(\bar{\ell}/\sqrt{n})}\bigl|\bar {a}^n\bigl(\widehat{X}^n (s)
\bigr)\bigr|\,ds \biggr]
\end{eqnarray*}
as required.
\end{pf*}

Toward the proof of Lemma~\ref{lem:minorization} we first prove that
$\hF_L^n(x)=L_n^{-1}\hF^n(L_n x)$ inherits the Lipschitz continuity of
$\hF^n$.
%
\begin{lem}%
\label{lem:PSD} There exists a finite positive constant $K$ (not depending
on $n$) such that
\[
\bigl|\hF_L^n(x)-\hF_L^n(y)\bigr|\leq
K|x-y|,\qquad x,y\in\bbR^d.
\]
\end{lem}
\begin{pf}Since, for each $n$,
$\bar{a}^n(0)$ is symmetric positive definite as is $\bar{a}$, these
matrices have strictly positive eigenvalues; see, for example, \cite{HoJ:94},
Theorem~7.2.1. Also, the eigenvalues of the square-root
matrix $L_n$ are the square roots of the eigenvalues of $\bar{a}^n(0)$.
Since $\bar{a}^n(0)\rightarrow\bar{a}$, the eigenvalues of $L_n$,
$(\lambda_1^n,\ldots,\lambda_d^n)$, converge to those of $L$,
$(\lambda
_1,\ldots,\lambda_d)$. The eigenvalues of the inverses $L_n^{-1}$ and
$L^{-1}$ are given by the reciprocals and, in turn, satisfy $(1/\lambda
_1^n,\ldots,1/\lambda_d^n)\rightarrow(1/\lambda_1,\ldots,1/\lambda_d)$.
In particular $|\!|\!|L_n|\!|\!|_2\rightarrow\|\!|L|\!\|_2$ and $\|\!|L_n^{-1}|\!\|
_2\rightarrow\|\!|L^{-1}|\!\|_2$ (where, following common notation, $|\!|\!|A|\!|\!|_2$ is the spectral norm of $A$; see \cite{HoJ:94}, Section~5.1.
Since the matrices are symmetric this norm is equal to the spectral
radius of the matrix, that is, to its maximal eigenvalue). By
definition of the matrix norm it then holds that
%
\renewcommand{\theequation}{\arabic{equation}}
\setcounter{equation}{55}
\begin{equation}
|L_n x-L_n y|\leq \|\!|L_n|\!
\|_2|x-y|\leq C_1\|\!|L|\!\|_2|x-y|,\qquad x,y\in
\bbR^d \label
{eq:interim20}
\end{equation}
for some finite positive constant $C_1$ where the
last inequality follows from the fact $\|\!|L_n|\!\|_2\rightarrow\|\!|L|\!\|
_2$ argued above. Similarly,
%
\begin{equation}
\bigl|L_n^{-1} x-L_n^{-1} y\bigr|\leq
C_2\bigl|\!\bigl|\! \bigl|L^{-1}\bigr|\!\bigr|\!\bigr|_2|x-y|,\qquad x,y\in
\bbR^d\label{eq:interim21}
\end{equation}
for a finite
positive constant $C_2$. Finally, using \eqref{eq:Flip} we have that
\begin{eqnarray*}
\bigl|L_n^{-1}\hF^n(L_n
x)-L_n^{-1}\hF^n(L_n y)\bigr|& \leq&\bigl|\!\bigl|\! \bigl|
L_n^{-1}\bigr|\!\bigr|\! \bigr|\bigl|\hF^n(L_n x)-
\hF^n(L_n y)\bigr|
\\
&\leq& C_2K_FC_1\bigl|\!\bigl|\! \bigl|L^{-1}\bigr|\!\bigr|\! \bigr|
_2\|\! |L|\!\|_2|x-y|,
\end{eqnarray*}
which completes the proof.
\end{pf}

\begin{pf*}{Proof of Lemma \protect\ref{lem:minorization}} We
consider first the
chain $\Phi^n$. Fix $K$ and let $\mathcal{K}=B_0(K)$. Let $\bar
{\mathcal
{K}}^n=L_n^{-1}\mathcal{K}$. By \eqref{eq:interim21}, there exists a
constant $\widetilde{K}$ not depending on $n$ such that
%
\begin{equation}
|x-y|\leq \widetilde{K}, \qquad x,y\in\bar{\mathcal{K}}^n.\label{eq:bark1}
\end{equation}

By Lemma~\ref{lem:PSD} there exist $\epsilon$ and $\delta$ not
depending on $n$ such that $|\hF_L^n(x)|\leq\epsilon+\delta|x-y|$ for
all $x\in\bbR^d$ and $y\in\bar{\mathcal{K}}^n$. Also, since $\hF
_L^n(0)=L_n^{-1}\hF^n(L_n 0)=0$ it satisfies also a linear growth
condition uniformly in $n$. Using \cite{qian2004representation}, Theorem~3.1 and \eqref{eq:bark1} we have that
\[
p(x,1,y)\geq\check{\epsilon},\qquad x,y\in\bar{\mathcal{K}}^n
\]
for some $\check{\epsilon}>0$ where $p(x,t,y)$ is the transition
density of $Z_L^n$ from $x$ to $y$ in time~$t$. In particular,
\[
\Pd_{\Phi^n}(x,\mathcal{E})\geq\int_{y\in\mathcal{E}\cap\bar
{\mathcal{K}}^n}p(x,1,y)\,dy
\geq\widetilde{\epsilon} \lambda\bigl(\bar {\mathcal{K}}^n\bigr)
\mathcal{Q}^n(\mathcal{E}),
\]
where $\lambda$ is here the Lebesgue measure and
\[
\mathcal{Q}^n(\cdot)= \frac{\lambda(\cdot\cap\bar{\mathcal
{K}}^n)}{\lambda(\bar{\mathcal{K}}^n)}.
\]

Using the invariance of Lebesgue measure under invertible linear
transformations we have for any $\mathcal{E}\in\mathcal{B}(\bbR^d)$ that
\[
\mathcal{Q}^n\bigl(L_n^{-1}\mathcal{E}\bigr)=
\frac{\lambda(L_n^{-1}\mathcal
{E}\cap L_n^{-1}\mathcal{K})}{\lambda(L_n^{-1}\mathcal{K})}=\frac
{\operatorname{det}(L_n^{-1})\lambda(\mathcal{E}\cap\mathcal
{K})}{\operatorname{det}(L_n^{-1})\lambda(\mathcal{K})},
\]
where $\operatorname{det}(L_n^{-1})>0$ is here the determinant of the positive
definite matrix $L_n^{-1}$, and we use the simple fact that
$(L_n^{-1}\mathcal{E})\cap(L_n^{-1} \mathcal
{K})=L_n^{-1}(\mathcal
{E}\cap\mathcal{K})$. Since $L_n\rightarrow L$, it also holds that
$\operatorname{det}(L_n^{-1})=(\operatorname{det}(L_n))^{-1}\rightarrow(\operatorname{det}(L))^{-1}=\operatorname{det}(L^{-1})>0$
so that there exists $\varepsilon>0$ (not depending on $n$) such that
\[
\lambda\bigl(\bar{\mathcal{K}}^n\bigr)=\operatorname{det}\bigl(L_n^{-1}
\bigr)\lambda(\mathcal{K})\geq \varepsilon.
\]
Let $\epsilon=\widetilde{\epsilon}\varepsilon$. Defining the measure
\[
\mathcal{Q}(\cdot)=\frac{\lambda(\cdot\cap\mathcal
{K})}{\lambda
(\mathcal{K})},
\]
we conclude that
\[
\Pd_{\Phi^n}(x,\mathcal{E})\geq\widetilde{\epsilon} \lambda\bigl(\bar {
\mathcal{K}}^n\bigr)\mathcal{Q}^n(\mathcal{E})=\epsilon
\mathcal {Q}(L_n\mathcal{E}),\qquad x\in\bar{\mathcal{K}}^n,
\mathcal{E}\in \mathcal {B}\bigl(\bbR^d\bigr).
\]

The result for $\Pd_{\psi^n}$ follows immediately from the above. Indeed,
\[
\Pd_{\psi^n}(x,\mathcal{E})=\Pd_{\Phi
^n}\bigl(L_n^{-1}x,L_n^{-1}
\mathcal {E}\bigr)\geq\epsilon\mathcal{Q}(\mathcal{E}), \qquad x\in\mathcal {K}, \mathcal
{E}\in\mathcal{B}\bigl(\bbR^d\bigr),
\]
which completes the proof.
\end{pf*}

\begin{pf*}{Proof of Lemma \protect\ref{lem:fromdisctocont}} This
argument is
almost identical to the proof in \cite{galtchouk2010geometric}, page
27. Under condition \eqref{eq:expobound1},
Dynkin's formula holds up to $t$, that is,
\[
\Ex_y\bigl[V\bigl(\hY^n(t)\bigr)\bigr]= V(y)+
\Ex_y \biggl[\int_0^t \mathcal
{A}^n V\bigl(\hY^n(s)\bigr)\,ds \biggr];
\]
see, for example, \cite{klebaner2005introduction}, Theorem~6.3. Setting
\[
g(t)=\Ex_y\bigl[V\bigl(\hY^n(t)\bigr)\bigr]\quad\mbox{and}\quad
 h(t)=\Ex_y\bigl[\mathcal{A}^nV\bigl(\hY
^n(t)\bigr)\bigr]+\delta g(t),
\]
we have that $h(t)\leq b\mathbh{1}_{\overline{B}_0(K)}(y)$ ($b$ and
$\delta$ as
in Assumption~\ref{asum:L}) and
\[
\dot{g}(t)=-\delta g(t)+h(t).
\]
Solving this differential equation we get
\begin{eqnarray*}
g(t)& =&g(0)e^{-\delta t}+\int_0^t
e^{\delta
(t-s)}h(s)\,ds\leq g(0)e^{-\delta t}+b\mathbh{1}_{\overline
{B}_0(K)}(y)
\frac
{1-e^{-\delta}}{\delta}
\\
& =& V(y)e^{-\delta t}+b\mathbh{1}_{\overline{B}_0(K)}(y)\frac
{1-e^{-\delta
}}{\delta}.
\end{eqnarray*}
Setting $\gamma=1- e^{-\delta}$ and $\bar{b}=b\frac{1-e^{-\delta
}}{\delta}$ we have the statement of the lemma.
\end{pf*}
\end{appendix}

\section*{Acknowledgments} The author is grateful to Junfei Huang and
to an anonymous referee for their careful reading of this paper and for
their numerous insightful comments.

%



\printaddresses

\end{document}